\def\E{\end{document}}
\documentclass[12pt]{article}
\usepackage{amssymb,amsmath,color}
\usepackage{enumitem}
\usepackage{url}
\allowdisplaybreaks[4]

\def\pf{\it{Proof.}\rm\quad}

\newtheorem{defn}{Definition}[section]
\newtheorem{thm}{Theorem}[section]
\newtheorem{ex}{Example}[section]
\newtheorem{cor}{Corollary}[section]

\newtheorem{lem}{Lemma}[section]
\newtheorem{rem}{Remark}[section]

\newtheorem{proposition}{Proposition}[section]

\catcode`@=11 \@addtoreset{equation}{section} \catcode`@=12

\begin{document}
\title{Sphere version of Banach--Kadec--Paley theorem}
\author{Qingjin Cheng$^a$, Wuyi He$^{b}$, Bo Xiang$^{a}$\thanks{Emails: qjcheng@xmu.edu.cn (Q. J. Cheng), hewy@cqupt.edu.cn (W. Y. He), xiangbo7501@gmail.com (B. Xiang)}\\
\small $^a$ School of Mathematical Sciences, Xiamen University, Xiamen 361005, China\\
\small $^b$ School of Mathematics and Statistics, Chongqing University of Posts and  Telecommunications, \\ \small Chongqing  400065, China}
\date{}
\maketitle 
\begin{abstract}
We study almost Lipschitz embeddings in the quantitative geometry of unit spheres of classical \(L_p\)-spaces.
Our main result is a sphere version of the classical
Banach--Kadec--Paley theorem: for \(1\le p,q<\infty\), the unit sphere
\(S_{L_q}\) admits an almost Lipschitz embedding into \(L_p\) if and
only if \(L_q\) admits a linear embedding into \(L_p\).
Consequently, the almost Lipschitz geometry of \(L_p\)-spheres
completely determines the underlying linear structure of 
\(L_p\)-spaces.

The proof relies on quantitative estimates for moduli of continuity of
sphere mappings. These estimates are derived from
Kalton--Randrianarivony's concentration inequalities,
Mendel--Naor's metric cotype inequality, and Naor's sharp metric
\(X_p\) inequality. They also yield several sharp results on the Hölder
geometry of \(L_p\)-spheres.
\end{abstract}
\maketitle
\tableofcontents

\section{Introduction}

\quad

A Banach space is, by its nature, also a metric space. A fundamental question of nonlinear geometry is to determine to what extent the metric geometry of a Banach space already determines its
linear structure.

Since Ribe's rigidity theorem~\cite{Ribe1976}, this question has become one
of the central themes of nonlinear Banach space theory.  Over the past five
decades it has led to a rich theory based on uniform, Lipschitz and coarse
embeddings between Banach spaces; see, for instance,
\cite{BL2000,AK2016,BL2026}.  One of the principal objectives of the
Ribe program is to identify metric invariants which recover linear properties
of Banach spaces.

A prototypical case is provided by the classical $L_p$ spaces, that is, $L_p=L_p([0,1],\lambda)$ where $\lambda$ is the Lebesgue measure.  One of the
oldest problems in the geometry of Banach spaces asks for which
$1\le p,q<\infty$ the space $L_q$ embeds linearly into $L_p$.  Chapter 12 of Banach's book
\emph{Th\'eorie des op\'erations lin\'eaires} (1932) \cite{Banach1932}
is devoted to this problem.  Banach proved that if $L_q$ is linearly isomorphic to a
subspace of $L_p$, then necessarily
\[
p\le q\le 2
\qquad\text{or}\qquad
2\le q\le p,
\]
and also showed that $L_2$ embeds linearly into every $L_p$.  He asked whether
the two remaining cases
\[
p<q<2
\qquad\text{and}\qquad
2<q<p
\]
could occur.  The first question was answered affirmatively by Kadec
\cite{Kadec1958}, who proved that every $L_q$ with $1\le p<q<2$ embeds
linearly into $L_p$.  The second was settled negatively by Paley
\cite{Paley1936}, who proved that $L_q$ cannot be linearly embedded into
$L_p$ whenever $2<q<p$.  These results culminated in the celebrated
Banach--Kadec--Paley classification theorem; see, for example,
\cite[Chapter~6]{AK2016}.  It states that, for $1\leq p\neq q<\infty$,
\[
L_q ~ {\rm embeds~ linearly ~into}~ L_p\Longleftrightarrow\quad q=2 \quad {\rm or} \quad 1\le p\le q\le2.
\]
Thus the linear embeddability relation among the classical $L_p$ spaces is
completely determined.

A natural question is whether this classification remains valid beyond the
linear category.  This question has played a decisive role in the development
of nonlinear Banach space geometry.  For uniform embeddings, however, the situation changes
completely: uniform geometry no longer determines the underlying linear
structure. For Lipschitz embeddings, differentiability methods originating in the work of
Mankiewicz and further developed by Johnson, Lindenstrauss, Preiss and
Schechtman reduce the classification problem to the linear theory; see
\cite{Lindenstrauss1964,JLPS2002}.  Hence the Lipschitz classification of the
classical $L_p$ spaces essentially coincides with the Banach--Kadec--Paley
classification.

The differentiability approach, however, does not provide a purely metric
explanation of the classification.  The search for such explanations has had a
major impact on the Ribe program.  Indeed, each forbidden range in the
nonlinear classification of $L_p$ spaces has led to the discovery of a
fundamental metric invariant.  The range
\[
q<\min\{p,2\}
\]
is related to the development of metric type, while the range
\[
q>\max\{p,2\}
\]
is explained by metric cotype, in particular by the theory of Mendel and
Naor~\cite{MN2008}.  The remaining range
\[
2<q<p
\]
was addressed by Naor and Schechtman \cite{NS2016} through their metric \(X_p\) inequalities. Naor \cite{Naor2016} subsequently established the sharp form of this inequality with the optimal scaling parameter, yielding sharp restrictions on snowflake embeddings. 
Thus the
study of nonlinear embeddings between the classical $L_p$ spaces has not only
determined the linear classification in metric terms, but has also produced
some of the most powerful tools in modern nonlinear Banach space geometry.

The present paper investigates the corresponding problem in the setting of unit
spheres.  More precisely, we ask whether a Banach-space classification theorem
can be determined  from nonlinear embeddings defined only on unit spheres.  This
question is substantially different from the whole space problem.  Unit
spheres carry no linear structure, and tools such as differentiation,
linearization and ultrapower arguments are no longer directly available.
Moreover, the classical uniform theory of spheres shows that the naive analogue
of the above classification is false.  Mazur's sphere theorem asserts that the
unit spheres of all infinite-dimensional separable $L_p$ spaces are uniformly
homeomorphic~\cite{Mazur1929}.  In fact, the Mazur map provides explicit
bi-H\"older homeomorphisms between these spheres.  Consequently, the uniform and Hölder geometries of Banach spheres are too
coarse to characterize the linear embeddability of the classical $L_p$
spaces.
This naturally leads to the search for a nonlinear category in which a
sphere analogue of the Banach--Kadec--Paley classification theorem
remains valid.

The main result of the present paper shows that almost Lipschitz
embeddings provide precisely the appropriate nonlinear category for such
a classification. Roughly speaking, almost Lipschitz embeddings are
arbitrarily close to being Lipschitz, up to logarithmic corrections, and
therefore constitute a natural intermediate category between Lipschitz
and uniform embeddings.

\begin{thm}[Sphere version of  the Banach--Kadec--Paley theorem]
\label{thm:intro-spherical-BKP}
Let $1\le p,q<\infty$. Then $S_{L_q}$ embeds  almost  Lipschitzly into $L_p$ if and only if $L_q$ embeds linearly into $L_p.$
\end{thm}

Thus, the almost Lipschitz geometry of unit spheres completely
determines the linear structure of classical \(L_p\)-spaces.
Accordingly, Theorem~\ref{thm:intro-spherical-BKP} may be viewed as a
sphere version of the classical Banach--Kadec--Paley theorem.

The proof is based on a quantitative theory of sphere mappings. The
central idea is to regard the moduli of continuity
\((\omega_f,\rho_f)\) as quantitative invariants of sphere embeddings
and to derive sharp restrictions on their possible behavior.

Our approach transfers several fundamental tools from nonlinear Banach
space geometry to the setting of unit spheres, including the
Kalton--Randrianarivony concentration inequality~\cite{KR2008},
Mendel--Naor metric cotype inequality~\cite{MN2008}, and Naor's sharp
metric \(X_p\) inequality~\cite{Naor2016,NS2016}. These yield the
quantitative estimates that form the technical core of the paper.

As applications, we establish sharp results on the Hölder geometry of
\(L_p\)-spheres, including a sphere version of Naor's bi-\(\theta\)-Hölder
embedding theorem, and develop the almost Lipschitz geometry of unit
spheres, culminating in Theorem~\ref{thm:spherical-KPP}.

The emphasis of the present work is not merely on obtaining quantitative
estimates for sphere mappings, but on showing that such estimates lead to
classification results for Banach spheres. We hope that this viewpoint
provides a first step towards a quantitative theory of Banach spheres and,
more broadly, towards a sphere version of the Ribe program.

The paper is organized as follows.  In Section~2 we recall basic notions on
moduli of continuity, Lipschitz embeddings, H\"older embeddings.
Section~3 develops estimates for moduli of continuity of sphere mappings,
derived from concentration inequalities, metric cotype and metric $X_p$
inequalities.  Section~4 applies these estimates to the H\"older geometry of
classical $L_p$ spheres, including sharp estimates for the Mazur map, optimal Hölder exponents for sphere embeddings, and a
sphere version of Naor's bi-$\theta$-H\"older embedding theorem.  Section~5 is devoted to almost Lipschitz embeddings and culminates in
the proof of a sphere version of the Banach--Kadec--Paley theorem,
showing that, for \(L_p\)-spaces, the almost Lipschitz geometry of unit
spheres completely determines the underlying linear structure. The
section concludes with a pair of almost Lipschitz homeomorphic Banach
spaces that are not Lipschitz homeomorphic, demonstrating that such
rigidity is peculiar to \(L_p\)-spaces rather than a general phenomenon.

\section{Basic notions and nonlinear embeddings}

\quad
In this paper, our results apply equally well to any infinite-dimensional Lebesgue
function space $L_p(\mu)$, but for concreteness we fix (as usual) the space $L_p$ to be equal to $L_p([0,1],\lambda)$, where $\lambda$ is the Lebesgue measure. Banach spaces are assumed to be over real scalars unless stated otherwise, though our results hold
true mutatis mutandis for complex Banach spaces as well.

We begin by recalling several standard notions from nonlinear Banach space
theory.  Let $(\mathcal{M},d_\mathcal{M})$ and $(\mathcal{N},d_\mathcal{N})$ be two metric spaces. Given a map  $f : \mathcal{M}\rightarrow \mathcal{N}$, define its extension modulus $\omega_f:[0,\infty)\rightarrow [0,\infty]$ by setting
\[\omega_{f}(t)=\sup\{d_\mathcal{N}(f(x),f(y)) : ~d_\mathcal{M}(x,y)\leq t\}\]
and its compression modulus \[\rho_f(t):=\inf\{d_\mathcal{N}(f(x),f(y)): ~d_\mathcal{M}(x,y)\geq t\}\]
respectively, where we adopt the conventions
\[
\sup(\emptyset)=0,
\qquad
\inf(\emptyset)=+\infty.
\]
Consequently,
\[
\rho_f\!\left(d_{\mathcal M}(x,y)\right)
\le
d_{\mathcal N}(f(x),f(y))
\le
\omega_f\!\left(d_{\mathcal M}(x,y)\right)
\]
for all \(x,y\in\mathcal M\).

\medskip

(i)
The map \(f\) is said to be \emph{uniformly continuous} if
\[
\lim_{t\to0^+}\omega_f(t)=0.
\]
It is called a \emph{uniform embedding} if, in addition,
\(\rho_f(t)>0\) for every \(t>0\).
A bijective uniform embedding is called a
\emph{uniform homeomorphism}.
The metric spaces \(\mathcal M\) and \(\mathcal N\) are said to be
\emph{uniformly homeomorphic} if there exists a uniform homeomorphism
between them.

\medskip

(ii)
The map \(f\) is called \emph{Lipschitz} if
\(\omega_f(t)\le Ct\) for some constant \(C>0\).
The least such constant,
\[
{\rm Lip}(f)
=
\sup_{x\neq y}
\frac{d_{\mathcal N}(f(x),f(y))}
{d_{\mathcal M}(x,y)},
\]
is called the \emph{Lipschitz constant} of \(f\).

\medskip

(iii)
Let \(\alpha>0\).
The map \(f\) is called
\emph{\(\alpha\)-H\"older} if
\(\omega_f(t)\le Ct^\alpha\) for some constant \(C>0\).
The least admissible constant,
\[
[f]_\alpha
=
\sup_{x\neq y}
\frac{d_{\mathcal N}(f(x),f(y))}
{d_{\mathcal M}(x,y)^\alpha},
\]
is called the \emph{\(\alpha\)-H\"older constant} of \(f\).

Given \(\alpha,\beta>0\), a bijection
\(f:\mathcal M\to\mathcal N\)
is called an
\emph{\((\alpha,\beta)\)-H\"older homeomorphism}
if \(f\) is \(\alpha\)-H\"older and
\(f^{-1}\) is \(\beta\)-H\"older.

\medskip

Throughout the paper, we use the following notation for embeddings:
\[
\begin{aligned}
X\hookrightarrow Y
&\quad &&\text{linear embedding},\\
X\stackrel{\rm Lip}{\hookrightarrow}Y
&\quad &&\text{Lipschitz embedding},\\
X\stackrel{\rm AL}{\hookrightarrow}Y
&\quad &&\text{almost Lipschitz embedding},\\
X\stackrel{\rm u}{\hookrightarrow}Y
&\quad &&\text{uniform embedding}.
\end{aligned}
\]

We also adopt the following conventions.
The symbols $C,\ C_1,\ C_2,\ \ldots,\ C_p,\ C_{p,q},\ C_{p,q,\theta},\ \ldots$ denote positive constants whose values may vary from line to line.
The dependence of a constant on parameters is indicated by its
subscripts whenever convenient.

For two nonnegative quantities \(A\) and \(B\), we write $A\lesssim B$
if there exists a constant \(C>0\) such that \(A\le CB\), and write $A\gtrsim B$
if \(B\lesssim A\).
We write $A\asymp B$ if both $A\lesssim B$ and $A\gtrsim B$
hold.

\section{Estimates for Moduli of Continuity of Sphere Mappings }

\quad

\subsection{A Lower Bound for the Modulus of Continuity}

\quad

Throughout the paper, the unit sphere $S_X$ is regarded as a metric
subspace of the Banach space $X$, equipped with the induced metric
\[
d(x,y)=\|x-y\|_X.
\]

The purpose of this section is to establish several estimates for moduli of
continuity on unit spheres. We establish a linear lower bound for the modulus of
continuity of injective uniformly continuous maps on unit spheres, which
will be used in the sequel.

\begin{lem}\label{HQ}
Let $X$ and $Y$ be infinite-dimensional Banach spaces, and let
$f:S_X\to Y$ be uniformly continuous. Then its modulus of continuity
$\omega_f$ satisfies
\begin{equation}\label{QH00}
\omega_f(s+t)
\le
2\bigl(\omega_f(s)+\omega_f(t)\bigr)
\end{equation}
for all $s,t\ge0$ with $s+t\le2$.
\end{lem}
\pf
Since any two points of $S_X$ can be joined by a
spherical arc of length at most $4$
(see \cite[3.5. Theorem]{Schaffer1967}), uniform continuity of $f$ implies
that $f(S_X)$ is bounded; hence $\omega_f(2)<\infty$.
Fix $s,t\ge0$ with $s+t\le2$, and let $\varepsilon>0$.
Choose $x_0,y_0\in S_X$ such that
\[
\|x_0-y_0\|_X\le s+t
\]
and
\[
\omega_f(s+t)-\varepsilon
\le
\|f(x_0)-f(y_0)\|_Y.
\]

Choose a two-dimensional subspace \(E\subseteq X\) containing \(x_0\) and
\(y_0\). More precisely, if \(\operatorname{span}\{x_0,y_0\}\) is
two-dimensional, let
\(
E=\operatorname{span}\{x_0,y_0\}.
\)
If \(y_0=\pm x_0\), choose \(z_{0}\notin \operatorname{span}\{x_0\}\) and put
\(
E=\operatorname{span}\{x_0,z_{0}\}.
\)
Since \(E\) is endowed with the restriction of the norm of \(X\), we have
\(
S_E=E\cap S_X.
\)
In particular, \(x_0,y_0\in S_E\).
Denote by $\overset{\frown}{xy}$ the shortest arc of $S_E$ joining
$x$ and $y$, and let $l(x,y)$ be its arc length.
By \cite[3.5. Theorem]{Schaffer1967} again,
\[
\|x-y\|_X
\le
l(x,y)
\le
2\|x-y\|_X,
\qquad x,y\in S_E.
\]
Hence
\[
l(x_0,y_0)\le 2(s+t).
\]

Since $(S_E,l)$ is metrically convex, there exists
$z\in S_E$ such that
\[
l(x_0,z)
=
\frac{s}{s+t}\,l(x_0,y_0),
\qquad
l(z,y_0)
=
\frac{t}{s+t}\,l(x_0,y_0).
\]
Therefore,
\[
l(x_0,z)\le 2s,
\qquad
l(z,y_0)\le 2t.
\]
Consequently,
\[
\|x_0-z\|_X\le 2s,
\qquad
\|z-y_0\|_X\le 2t.
\]
Thus
\[
\begin{aligned}
\|f(x_0)-f(y_0)\|_Y
&\le
\|f(x_0)-f(z)\|_Y
+
\|f(z)-f(y_0)\|_Y\\
&\le
\omega_f(2s)
+
\omega_f(2t).
\end{aligned}
\]
Since $\varepsilon>0$ is arbitrary,
\begin{equation}\label{QW1}
\omega_f(s+t)
\le
\omega_f(2s)+\omega_f(2t).
\end{equation}

More generally, the same arc-partition argument gives

\begin{equation} \label{QW00} \omega_f\left(\sum_{i=1}^n s_i\right)
\le
\sum_{i=1}^n\omega_f(2s_i) 
\end{equation}
for every $n\in\mathbb N$ and all $s_1,\ldots,s_n\ge0$
with $s:=\sum_{i=1}^n s_i\le2$.
The case $s=0$ is immediate.
If $s>0$, take distinct $x,y\in S_X$ with
$\|x-y\|_X\le s$. Choose a shortest arc $\gamma$
joining $x$ to $y$ in the unit sphere of a
two-dimensional subspace containing them, and denote
its length by $L$. As above,
\[
L\le 2\|x-y\|_X\le 2s.
\]
Choose successive points
\[
x=z_0,z_1,\ldots,z_n=y
\]
on $\gamma$ so that the $i$-th subarc has length
$Ls_i/s$, allowing repeated points when $s_i=0$.
 Then
\(
\|z_i-z_{i-1}\|_X\le\frac{Ls_i}{s}\le2s_i,
\)
and hence
\[
\|f(x)-f(y)\|_Y
\le
\sum_{i=1}^n\|f(z_i)-f(z_{i-1})\|_Y
\le
\sum_{i=1}^n\omega_f(2s_i).
\]
The same estimate is trivial when $x=y$.
Taking the supremum over all such pairs proves~(\ref{QW00}).

Applying \eqref{QW00} with
\[
s_1=s_2=\frac{s}{2},
\qquad
s_3=s_4=\frac{t}{2},
\]
we obtain
\[
\begin{aligned}
\omega_f(s+t)
&=
\omega_f\Bigl(
\frac{s}{2}+\frac{s}{2}
+
\frac{t}{2}+\frac{t}{2}
\Bigr)\\
&\le
\omega_f(s)+\omega_f(s)
+
\omega_f(t)+\omega_f(t)\\
&=
2\bigl(\omega_f(s)+\omega_f(t)\bigr),
\end{aligned}
\]
which proves \eqref{QH00}.\hfill$\square$

\begin{cor}\label{Qcor}
Let $X$ and $Y$ be infinite-dimensional Banach spaces, and let $f:S_X\to Y$ be uniformly continuous. Then for every $n\in \mathbb{N}$ and $s\ge0$ with $ns\leq 2$, 
\[
\omega_f(ns)\le 2 n\,\omega_f(s).
\]
\end{cor}

\pf
Let $n\in\mathbb N$ and $s\ge0$ with $ns\le2$.
Applying~\eqref{QW00} in Lemma~\ref{HQ} to $2n$ terms, each equal to $s/2$,
we obtain
\[
\omega_f(ns)
=
\omega_f\left(\sum_{i=1}^{2n}\frac{s}{2}\right)
\le
\sum_{i=1}^{2n}\omega_f(s)
=
2n\omega_f(s).
\]
\hfill$\square$

\begin{lem}\label{Q1}
Let $X$ and $Y$ be infinite-dimensional Banach spaces, and let $f:S_X\to Y$ be an injective uniformly continuous map. Then
\[
\omega_f(\varepsilon)\ge \frac{1}{4}\,\omega_f(1)\,\varepsilon
\quad \text{for all }\quad 0\le \varepsilon\le 2.
\]
\end{lem}
\pf
For $0<\varepsilon\le 2$, write $\omega_f(\varepsilon)=C(\varepsilon)\,\varepsilon$.
Choose $n\in\mathbb{N}$ such that
\[
\frac{1}{1+n}\le \frac{\varepsilon}{2} \le \frac{1}{n}.
\]
By the quasi-subadditivity of $\omega_f$ (Lemma~\ref{HQ}),
for any $m\in\mathbb{N}$ we have
\[
\omega_f(1)=\omega_f\Bigl(2m\cdot \frac{1}{2m}\Bigr)
\le 2m\,\omega_f\Bigl(\frac{1}{m}\Bigr).
\]
Taking $m=n$, it follows that
\[
\frac{1}{2n}\,\omega_f(1)\le \omega_f\Bigl(\frac{1}{n}\Bigr)
\le \omega_f(\varepsilon)=C(\varepsilon)\,\varepsilon
\le C(\varepsilon)\,\frac{2}{n}.
\]
Hence
\[
C(\varepsilon) \ge \frac{n}{2}\cdot \frac{\omega_f(1)}{2n} = \frac{1}{4}\,\omega_f(1),
\]
which proves the claim.\hfill$\square$

\subsection{Estimates via Concentration Inequalities }

\quad

We begin this section by introducing an alternative graph technique of the Gorelik principle.
The following  so called \emph{Kalton-Randrianarivony's graphs}  have
been introduced and developed by Kalton and Randrianarivony in \cite{KR2008}, which are also called \emph{Hamming graphs} in \cite{Fovelle2022, Fovelle2025}.\par
Let $\mathbb{M}$ be an infinite subset of $\mathbb{N}$. We denote by
$[\mathbb{M}]^{\omega}$ the set of infinite subsets of $\mathbb{M}$. For $M\in [\mathbb{N}]^{\omega}$ and $k\in\mathbb{N}$, let

$$[\mathbb{M}]^k=\Big\{\overline{n}=(n_1,\cdots,n_k)\in \mathbb{M}^{k},~ n_1<n_2<\cdots<n_k\Big\},$$
Then we equip $[\mathbb{M}]^k$ with
the Hamming distance $d_\mathbf{H}$, i.e.,
\begin{equation*}
  d_\mathbf{H}(\overline{n},\overline{m})=\left|\big\{j: n_j\neq m_j\big\}\right|.
\end{equation*}
for all $\overline{n}=(n_1,\cdots,n_k),\overline{m}=(m_1,\cdots,m_k)\in  [\mathbb{M}]^k$. Here $|A|$ denotes the cardinality of  a set $A\subset \mathbb{N}$.\par

\begin{defn} \cite{Fovelle2022, Fovelle2025} \label{Fovelle}
 Let $(\mathcal{M},d_\mathcal{M})$ be a metric space, $\lambda>0$, $p\in (1,\infty]$. We say that $\mathcal{M}$ has property \emph{$\lambda$-${\rm HFC}_{p}$ (Hamming Full Concentration)} if, for every $k\in\mathbb{N}$ and  every Lipschitz function $f : \left([\mathbb{N}]^{k}, d_\mathbf{H}\right)\to (\mathcal{M},d_\mathcal{M})$, one can find $\mathbb{M}\in[\mathbb{N}]^{\omega}$ such that
 $$d_\mathcal{M}\big(f(\overline{n}), f(\overline{m})\big)\leq\lambda k^{\frac{1}{p}} {\rm Lip} (f),~~\forall \, \overline{n},\overline{m}\in[\mathbb{M}]^{k}.$$
We say that $\mathcal{M}$ has property ${\rm HFC}_{p}$ if $\mathcal{M}$ has property $\lambda$-${\rm HFC}_{p}$ for some $\lambda>0$.
\end{defn}
\par
 The
 fundamental result of Kalton and Randrianarivony  (Theorem 6.1 in \cite{KR2008}) can be rephrased as
 follows.
\begin{thm}\cite{KR2008} \label{KR2008}
 Let $1<q<\infty$. Every reflexive $q$-{\rm AUS} Banach space has property ${\rm HFC}_{q}$.
\end{thm}

 According to \cite{JLPS2002}, we note that $\ell_{p}$ is $p$-AUS and $L_{p}$ is $\min\{p,2\}$-AUS (see also \cite[\text{ p. }114\text{--}117]{Milman1971}) for $1<p<\infty$.
It is important to mention here that  Kalton \cite[Remarks on p.170]{kalton2013} proved the existence of a Banach space $X$ that has property ${\rm HFC}_q$, but it is not $q$-AUSable for $1<q<\infty$.

 Let $X$ be a Banach space and $t>0$. For $x\in S_{X}$ and $E$ a closed linear subspace of $X$, we define
\begin{equation*}
  \overline{\rho}_X(t)=\sup_{x\in S_X}\inf_{\text{dim} (X/E)<\infty}\sup_{y\in S_E}\|x+ty\|-1.
\end{equation*}
$X$ is said to be \emph{asymptotically uniformly smooth} (in short AUS) if $\lim_{t\rightarrow 0}\frac{\overline{\rho}_X(t)}{t}=0$  . In particular, if there is $c>0$ and $q\in (1,\infty)$ such that $\overline{\rho}_X(t)\leq ct^q$ for $t>0$ then $X$ is called $q$-\emph{asymptotically uniformly smooth} (in short $q$-AUS).
 These notions were first introduced by Milman \cite{Milman1971} under different
notations and names. We refer to \cite{JLPS2002} for a survey about the relevant results concerning uniform asymptotic smoothness.



\begin{thm}\label{HFC}
Let \(1\le p<q<\infty\), and suppose that \(Y\) has property
\({\rm HFC}_q\). If
\[
f:S_{\ell_p}\longrightarrow Y
\]
is a uniform embedding, then there exists a constant \(K>0\) such that
\[
\omega_f(\varepsilon)\ge K\varepsilon^{p/q}
\]
for every \(0\le \varepsilon\le 2\).
\end{thm}
\pf
Fix \(0<\varepsilon\le2\). Choose \(k\in\mathbb \{2,3,\cdots\}\) such that
\[
\frac1{k}< \left(\frac{\varepsilon}{2}\right)^p\le\frac2{k} .
\]
Define
\[
\phi:[\mathbb N]^k\to S_{\ell_p},\qquad
\phi(n_1,\ldots,n_k)=k^{-1/p}\sum_{i=1}^k e_{n_i}.
\]
Then \(\|\phi(\bar n)\|_p=1\) for every
\(\bar n\in[\mathbb N]^k\). Let \(g=f\circ\phi\).

If \(d_{\mathbf H}(x,y)=1\), then
\[
\|\phi(x)-\phi(y)\|_p=2^{1/p}k^{-1/p},
\]
and hence
\[
\|g(x)-g(y)\|
\le
\omega_f(2^{1/p}k^{-1/p}).
\]

Note that for any
\(\bar n,\bar m\in[\mathbb N]^k\) with
\(d_{\mathbf H}(\bar n,\bar m)=r\), there exists a chain
\[
\bar n=x_0,x_1,\ldots,x_r=\bar m
\]
such that \(d_{\mathbf H}(x_{i-1},x_i)=1\) for every
\(i=1,\ldots,r\). Therefore,
\[
\begin{aligned}
\|g(\bar n)-g(\bar m)\|
&\le
\sum_{i=1}^{r}\|g(x_i)-g(x_{i-1})\|  \\
&\le
r\,\omega_f(2^{1/p}k^{-1/p})\\
&=
d_{\mathbf H}(\bar n,\bar m)
\omega_f(2^{1/p}k^{-1/p}),
\end{aligned}
\]
and consequently
\[
\operatorname{Lip}(g)
\le
\omega_f(2^{1/p}k^{-1/p}).
\]

Since \(Y\) has property \({\rm HFC}_q\), there exist
\(\lambda>0\) and an infinite set \(\mathbb M\subseteq\mathbb N\) such
that
\[
\|g(\bar n)-g(\bar m)\|
\le
\lambda k^{1/q}\operatorname{Lip}(g),
\qquad
\bar n,\bar m\in[\mathbb M]^k .
\]

Choose disjoint \(\bar n,\bar m\in[\mathbb M]^k\). Then
\(\|\phi(\bar n)-\phi(\bar m)\|_p=2^{1/p}\), and hence
\[
\rho_f(2^{1/p})
\le
\lambda k^{1/q}
\omega_f(2^{1/p}k^{-1/p}).
\]
Since \(k^{-1/p}\le\varepsilon/2\), we have
\(2^{1/p}k^{-1/p}\le2^{1/p-1}\varepsilon\le\varepsilon\).
By the monotonicity of \(\omega_f\), it follows that
\[
\rho_f(2^{1/p})
\le
\lambda k^{1/q}\omega_f(\varepsilon),
\]
and therefore
\[
\omega_f(\varepsilon)
\ge
\frac{\rho_f(2^{1/p})}{\lambda}k^{-1/q}.
\]

Moreover, the choice of \(k\) gives
\((\varepsilon/2)^p\le2/k\), and hence
\[
k^{-1/q}\ge
2^{-(p+1)/q}\varepsilon^{p/q}.
\]
Combining the above estimates and using
\(\rho_f(2^{1/p})>0\), we obtain
\[
\omega_f(\varepsilon)\ge K\varepsilon^{p/q},
\]
where 
\[
K=\frac{\rho_f(2^{1/p})}{\lambda}2^{-(p+1)/q}>0 .
\]
\hfill$\square$

The concentration argument yields lower estimates for the expansion
modulus. To obtain complementary upper estimates for the compression
modulus, one needs a different approach based on metric cotype.

\subsection{Estimates via Metric Cotype Inequalities}

We begin by fixing some notation and conventions. For \(m\in\mathbb N\),
let \(\mathbb Z_m\) denote the additive group of integers modulo \(m\).
We equip \(\mathbb Z_m^n\) with the standard \(\ell_\infty\)-Cayley graph:
two distinct points \(x,y\in\mathbb Z_m^n\) are adjacent if
\[
x-y\in\{-1,0,1\}^n.
\]
The corresponding shortest-path metric is denoted by
\(d_{\mathbb Z_m^n}\). Equivalently,
\[
(\mathbb Z_m^n,d_{\mathbb Z_m^n})
=
(\mathbb Z^n,\|\cdot\|_\infty)/(m\mathbb Z)^n.
\]

We denote by \(\mu\) the normalized counting measure on
\(\mathbb Z_m^n\), and by \(\sigma\) the normalized counting measure on
\(\{-1,0,1\}^n\). For each \(j\in\{1,\ldots,n\}\), let \(e_j\) denote
the \(j\)-th canonical basis vector of \(\mathbb Z_m^n\).

\begin{defn}[\cite{MN2008}]\label{Cotype}
Let \((\mathcal M,d_{\mathcal M})\) be a metric space, and let
\(q,\Gamma>0\). We say that \(\mathcal M\) has \emph{metric cotype \(q\)}
with constant \(\Gamma\) if, for every \(n\in\mathbb N\), there exists
an even integer \(m\in\mathbb N\) such that every map
\(F:\mathbb Z_m^n\to\mathcal M\) satisfies
\begin{equation}\label{MC}
\begin{split}
&\sum_{j=1}^n
\int_{\mathbb Z_m^n}
d_{\mathcal M}\left(
F\left(x+\frac m2e_j\right),F(x)
\right)^q\,d\mu(x)
\\
&\qquad\le
\Gamma^qm^q
\int_{\{-1,0,1\}^n}
\int_{\mathbb Z_m^n}
d_{\mathcal M}\bigl(F(x+\eta),F(x)\bigr)^q
\,d\mu(x)\,d\sigma(\eta).
\end{split}
\end{equation}
For \(n\in\mathbb N\) and \(\Gamma>0\), we denote by
\(m_q(\mathcal M,n,\Gamma)\) the smallest even integer \(m\) for which
\eqref{MC} holds for every \(F:\mathbb Z_m^n\to\mathcal M\). If no such
integer exists, we set
\[
m_q(\mathcal M,n,\Gamma)=\infty.
\]
\end{defn}

The constructions of maps from discrete tori into unit spheres and the associated arguments in the following theorem
and in Theorem 3.5 below are inspired by Naor's proof of
\cite[Lemma 5.2]{Naor2014}. 

\begin{thm}\label{Cotype1}
Let \(2\le q<p<\infty\), and let \(X\) be a Banach space of cotype \(q\)
and nontrivial type. If
\[
f:S_{\ell_p}\longrightarrow X
\]
is a uniform embedding, then there exists \(K>0\) such that
\[
\rho_f(\varepsilon)\le K\omega_f(\varepsilon^{p/q}),
\qquad 0\le\varepsilon\le2 .
\]
\end{thm}

\pf
Since \(X\) has nontrivial type, it is \(K\)-convex. Hence, by
\cite[Theorem 4.1]{MN2008}, there exist constants
\(\Gamma,K_1>0\) such that
\[
m_q(X,n,\Gamma)\le K_1n^{1/q}
\]
for every \(n\in\mathbb N\). Moreover, by
\cite[Lemma 2.3]{MN2008},
\[
m_q(X,n,\Gamma)\ge \frac{n^{1/q}}{\Gamma}.
\]

Fix \(0<\varepsilon\le2\), and choose \(n\in\mathbb N\) such that
\[
\frac1{n+1}<\left(\frac{\varepsilon}{2}\right)^p\le\frac1n .
\]
Let \(m=m_q(X,n,\Gamma)\). Then
\[
\frac{n^{1/q}}{\Gamma}\le m\le K_1n^{1/q}.
\]

Define
\[
\phi(x)=n^{-1/p}
(e^{2\pi ix_1/m},\ldots,e^{2\pi ix_n/m},0,\ldots),
\qquad x\in\mathbb Z_m^n .
\]
Then \(\phi(x)\in S_{\ell_p}\).
Set \(g=f\circ\phi\).

For every \(\eta\in\{-1,0,1\}^n\), using
\(|e^{is}-e^{it}|\le |s-t|\), we have
\[
\|\phi(x+\eta)-\phi(x)\|_p
\le
n^{-1/p}
\Big(\sum_{j=1}^n(2\pi|\eta_j|/m)^p\Big)^{1/p}
\le\frac{2\pi}{m}.
\]
Hence,
\[
\int_{\{-1,0,1\}^n}\int_{\mathbb Z_m^n}
\|g(x+\eta)-g(x)\|^q\,d\mu(x)d\sigma(\eta)
\le
\omega_f(2\pi/m)^q .
\]

On the other hand,
\[
\left\|
\phi\left(x+\frac m2e_j\right)-\phi(x)
\right\|_p
=
2n^{-1/p},
\]
and therefore
\[
\sum_{j=1}^n\int_{\mathbb Z_m^n}
\left\|
g\left(x+\frac m2e_j\right)-g(x)
\right\|^q\,d\mu(x)
\ge
n\rho_f(2n^{-1/p})^q .
\]
Applying the metric cotype inequality gives
\[
n\rho_f(2n^{-1/p})^q
\le
\Gamma^qm^q\omega_f(2\pi/m)^q ,
\]
and hence
\begin{equation}\label{TAG35}
\rho_f(2n^{-1/p})
\le
\Gamma mn^{-1/q}\omega_f(2\pi/m)
\le
\Gamma K_1\omega_f(2\pi\Gamma n^{-1/q}).
\end{equation}

The choice of \(n\) implies
\[
n^{-1/p}<\varepsilon\le2n^{-1/p},
\]
and consequently
\[
n^{-1/q}<\varepsilon^{p/q}.
\]
Since \(\rho_f\) is nondecreasing, (\ref{TAG35}) yields
\begin{equation}\label{TAG36}
\rho_f(\varepsilon)
\le
\Gamma K_1\omega_f(2\pi\Gamma\varepsilon^{p/q}).
\end{equation}

Choose \(\varepsilon_0>0\) so small that
\[
\lceil2\pi\Gamma\rceil\varepsilon_0^{p/q}\le2 .
\]
For \(0<\varepsilon\le\varepsilon_0\), by Corollary~\ref{Qcor}
and the monotonicity of \(\omega_f\), we have
\[
\omega_f(2\pi\Gamma\varepsilon^{p/q})
\le
2\lceil2\pi\Gamma\rceil
\omega_f(\varepsilon^{p/q})
\le
2(2\pi\Gamma+1)
\omega_f(\varepsilon^{p/q}).
\]

Combining this with (\ref{TAG36}), we obtain
\[
\rho_f(\varepsilon)
\le
2\Gamma K_1(2\pi\Gamma+1)
\omega_f(\varepsilon^{p/q}).
\]

For \(\varepsilon_0<\varepsilon\le2\), the same conclusion follows after
enlarging the constant, since
\(\rho_f(\varepsilon)\le\rho_f(2)\) and
\(\omega_f\) is positive on the compact interval
\([\varepsilon_0^{p/q},2^{p/q}]\).

Therefore, there exists \(K>0\) such that
\[
\rho_f(\varepsilon)
\le
K\omega_f(\varepsilon^{p/q}),
\qquad 0\le\varepsilon\le2 .
\]
\hfill$\square$

\medskip

\subsection{Estimates via Metric $X_p$ Inequalities }

\quad

In what follows, given $n \in \mathbb{N}$ we shall denote the set $\{1,\dots,n\}$ by $[n]$.
The coordinate basis of $\mathbb{R}^n$ will be denoted by $e_1,\dots,e_n$, and for a sign vector $\varepsilon = (\varepsilon_1,\dots,\varepsilon_n) \in \{-1,1\}^n$ and a subset $S \subseteq [n]$, we shall use the notation
\begin{equation}
\label{eq:epsilon_S}
\varepsilon_S \stackrel{\mathrm{def}}{=} \sum_{j \in S} \varepsilon_j e_j.
\end{equation}

\begin{defn}
(\cite{NS2016, Naor2016}) Fix $p \in (0,\infty)$.
Metric space $(\mathcal{M},d_{\mathcal{M}})$ is said to be an $X_p$ metric space if there exists $\mathfrak{X} \in (0,\infty)$ such that for every $n \in \mathbb{N}$ and $k \in [n]$, there exists $m \in \mathbb{N}$ such that, denoting as usual $\mathbb{Z}_{2m} = \mathbb{Z}/(2m\mathbb{Z})$, every function $f : \mathbb{Z}_{2m}^n \to \mathcal{M}$ satisfies the following distance inequality:

\begin{equation}\label{XPM}
 \begin{split}   
&\left( \frac{1}{\binom{n}{k}} \sum_{\substack{S \subseteq [n] \\ |S|=k}} \mathbb{E}\left[ d_\mathcal{M}\left(f(x + m\varepsilon_S), f(x)\right)^p \right] \right)^{\frac{1}{p}}\\
&\le \mathfrak{X} m \left( \frac{k}{n} \sum_{j=1}^n \mathbb{E}\left[ d_\mathcal{M}\left(f(x + e_j), f(x)\right)^p \right] + \left( \frac{k}{n} \right)^{\frac{p}{2}} \mathbb{E}\left[ d_\mathcal{M}\left(f(x + \varepsilon), f(x)\right)^p \right] \right)^{\frac{1}{p}}.
\end{split}
\end{equation}

The expectations in (\ref{XPM}) are with respect to $(x,\varepsilon) \in \mathbb{Z}_{2m}^n \times \{-1,1\}^n$ chosen uniformly at random.
\end{defn}

As established by Naor in \cite{Naor2016}, the metric $X_{p}$ inequality with sharp scaling parameter holds true in $L_{p}$ for every $p\in [2,\infty)$.
\begin{thm}\label{XP}
(\cite[Theorem 4]{Naor2016}) Suppose that $k,m,n \in \mathbb{N}$ satisfy $k \in [n]$ and $m \geqslant \sqrt{n/k}$. Suppose also that $p \in [2,\infty)$. Then every $f : \mathbb{Z}_{8m}^n \to L_p$ satisfies
\begin{eqnarray*}
&& \left( \frac{1}{\binom{n}{k}} \sum_{\substack{S \subseteq [n] \\ |S|=k}} \mathbb{E}\left[ \| f(x + 4m \varepsilon_S) - f(x) \|_{L_p}^p \right] \right)^{\frac{1}{p}} \\
& \lesssim_p& m \left( \frac{k}{n} \sum_{j=1}^n \mathbb{E}\left[ \| f(x + e_j) - f(x) \|_{L_p}^p \right] + \left( \frac{k}{n} \right)^{\frac{p}{2}} \mathbb{E}\left[ \| f(x + \varepsilon) - f(x) \|_{L_p}^p \right] \right)^{\frac{1}{p}},
\end{eqnarray*}
where the expectations are taken with respect to $(x,\varepsilon) \in \mathbb{Z}_{8m}^n \times \{-1,1\}^n$ chosen uniformly at random.
\end{thm}

For \(2<q<p<\infty\), the metric \(X_p\) inequalities of Naor and Schechtman \cite{NS2016} provide a purely metric proof that \(L_q\) does not admit a bi-Lipschitz embedding into \(L_p\). This can be viewed as a Lipschitz variant of Paley's theorem in 1936.
In this section we aim to establish estimates for the moduli of continuity  on the sphere of $L_q$.

\begin{thm}\label{sphere-modulus-ineq}
Let \(2<q<p<\infty\). If
\(f:S_{L_q}\to L_p\) is a uniform embedding, then there exists a
constant \(C_p>0\), depending only on \(p\), such that for every
\(n\in\mathbb N\), every \(k\in[n]\), and every
\(m\ge\sqrt{n/k}\),
\[
\rho_f\!\left(2\Bigl(\frac{k}{n}\Bigr)^{1/q}\right)
\le
C_pmk^{1/p}
\omega_f\!\left(\frac{\pi}{4m}n^{-1/q}\right)
+
C_pm\Bigl(\frac{k}{n}\Bigr)^{1/2}
\omega_f\!\left(\frac{\pi}{4m}\right).
\]
\end{thm}

\pf
For every \(x=(x_1,\ldots,x_n)\in\mathbb Z_{8m}^n\), define
\[
\phi:\mathbb Z_{8m}^n\longrightarrow S_{\ell_q}
\]
by
\[
\phi(x)=
\left(
n^{-1/q}e^{2\pi x_1 i/8m},
\ldots,
n^{-1/q}e^{2\pi x_n i/8m},
0,0,\ldots
\right)\in S_{\ell_q}\subseteq S_{L_q}.
\] Set \(g=f\circ\phi\).

Let \(S\subseteq[n]\) with \(|S|=k\), and let
\(\varepsilon\in\{-1,1\}^n\).
Since
\(
|e^{it}-e^{is}|
=2\bigl|\sin\frac{t-s}{2}\bigr|
\le |t-s|,
\)
we have
\[
\|\phi(x+4m\varepsilon_S)-\phi(x)\|_{\ell_q}
=
2\Bigl(\frac{k}{n}\Bigr)^{1/q},
\]
and hence
\[
\|g(x+4m\varepsilon_S)-g(x)\|_{L_p}
\ge
\rho_f\!\left(
2\Bigl(\frac{k}{n}\Bigr)^{1/q}
\right).
\]

Similarly,
\[
\|\phi(x+e_j)-\phi(x)\|_{\ell_q}
\le
\frac{\pi}{4m}n^{-1/q},
\qquad
\|\phi(x+\varepsilon)-\phi(x)\|_{\ell_q}
\le
\frac{\pi}{4m},
\]
so that
\[
\|g(x+e_j)-g(x)\|_{L_p}
\le
\omega_f\!\left(\frac{\pi}{4m}n^{-1/q}\right),
\qquad
\|g(x+\varepsilon)-g(x)\|_{L_p}
\le
\omega_f\!\left(\frac{\pi}{4m}\right).
\]

Applying Theorem~\ref{XP} to \(g\) yields
\[
\rho_f\!\left(
2\Bigl(\frac{k}{n}\Bigr)^{1/q}
\right)
\le
C_pm
\left(
k\,\omega_f\!\left(\frac{\pi}{4m}n^{-1/q}\right)^p
+
\Bigl(\frac{k}{n}\Bigr)^{p/2}
\omega_f\!\left(\frac{\pi}{4m}\right)^p
\right)^{1/p}.
\]
Finally, using
\(
(a+b)^{1/p}\le a^{1/p}+b^{1/p}
\)
for \(a,b\ge0\), we obtain
\[
\rho_f\!\left(
2\Bigl(\frac{k}{n}\Bigr)^{1/q}
\right)
\le
C_pmk^{1/p}
\omega_f\!\left(\frac{\pi}{4m}n^{-1/q}\right)
+
C_pm
\Bigl(\frac{k}{n}\Bigr)^{1/2}
\omega_f\!\left(\frac{\pi}{4m}\right),
\]
which completes the proof.
\hfill$\square$

\medskip

The estimates established in this section form a unified
framework for studying quantitative sphere geometry. In the next section,
we illustrate their strength by applying them to several fundamental
problems concerning classical \(L_p\)-spaces.

\section{Applications to the Hölder Geometry of  $L_p$ Spheres  }

In this section we apply the spherical estimates  developed in
Section~3 to several classical problems in the geometry of \(L_p\)-spaces.
These applications demonstrate that quantitative information on unit
spheres already determines remarkable nonlinear, and even linear,
properties of the ambient Banach spaces.

\subsection{Sharpness of the Mazur map estimates}

The Mazur map plays a fundamental role in the nonlinear geometry of Banach
spaces. In his celebrated 1929 paper~\cite{Mazur1929}, Mazur introduced this
map to prove that the unit spheres of all complex $L_p$ spaces,
$1\le p<\infty$, are uniformly homeomorphic. Later,
Benyamini--Lindenstrauss~\cite[Theorem~9.1]{BL2000} and
Naor~\cite[p.~29]{Naor2014} established the standard estimates for its
moduli of continuity. They may be stated as follows.

\begin{thm}[Mazur's sphere theorem]\label{th1.1}
Let $(\Omega,\mu)$ be a measure space and let
$1\le p,q<\infty$.
The Mazur map
\[
M_{p,q}:L_p(\Omega,\mu,\mathbb C)\longrightarrow
L_q(\Omega,\mu,\mathbb C)
\]
is defined by
\[
M_{p,q}(f)=|f|^{\frac pq}\operatorname{sign}(f).
\]
Then its restriction
\[
M_{p,q}:S_{L_p(\Omega,\mu,\mathbb C)}
\longrightarrow
S_{L_q(\Omega,\mu,\mathbb C)}
\]
is a uniform homeomorphism satisfying

\begin{enumerate}

\item[(i)]
If $1\le q<p<\infty$, then
\[
\omega_{M_{p,q}}(t)\le
\frac pq\,t,
\qquad
\rho_{M_{p,q}}(t)\ge
2^{1-\frac pq}t^{\frac pq},
\qquad
0\le t\le2.
\]

\item[(ii)]
If $1\le p\le q<\infty$, then
\[
\omega_{M_{p,q}}(t)\le
2^{1-\frac pq}t^{\frac pq},
\qquad
\rho_{M_{p,q}}(t)\ge
\frac pq\,t,
\qquad
0\le t\le2.
\]

\end{enumerate}
\end{thm}

The purpose of this subsection is to show that all the estimates in
Theorem~\ref{th1.1} are sharp. We begin with the finite-dimensional
setting. Throughout, for $n\ge2$ and $1\le p<\infty$,
$\ell_p^n$ denotes $\mathbb C^n$ equipped with the $\ell_p$ norm.

\begin{thm}\label{th1.2}
Let $n\ge2$ and let
$m_{p,q}:S_{\ell_p^n}\to S_{\ell_q^n}$
be the Mazur map
\[
m_{p,q}(x)
=
\big(
|x_1|^{\frac pq}\operatorname{sign}(x_1),
\dots,
|x_n|^{\frac pq}\operatorname{sign}(x_n)
\big).
\]
Then all the estimates in
Theorem~\ref{th1.1} are sharp.
\end{thm}
\pf
Since
$\ell_p^n=L_p(\{1,\dots,n\},\mu_c)$,
where $\mu_c$ denotes the counting measure,
Theorem~\ref{th1.1} immediately yields the estimates.
It remains to prove their optimality.

\medskip

\noindent
{\bf Case 1.}
Suppose $1\le p<q<\infty$.

\smallskip

\noindent
{\it Optimality of the Hölder exponent.}

For sufficiently large $m$, let
\[
x_m=
\Big(
(1-\tfrac1m)^{1/p},
0,\dots,0,
(\tfrac1m)^{1/p}
\Big),
\qquad
y_m=(1,0,\dots,0).
\]
Then $x_m,y_m\in S_{\ell_p^n}$ and
\[
\lim_{m\to\infty}
\frac{\|m_{p,q}(x_m)-m_{p,q}(y_m)\|_q}
     {\|x_m-y_m\|_p^{p/q}}
=
\begin{cases}
2^{-1/q}, & p=1<q,\\[2mm]
1,        & 1<p<q.
\end{cases}
\]
(The limit follows from the first-order expansion
$(1-x)^\gamma=1-\gamma x+o(x)$ as $x\to0^+$.)

Now let $\alpha>\frac pq$ and write
$\alpha=(1+\beta)\frac pq$
for some $\beta>0$.
Since
$\|x_m-y_m\|_p\to0$,
\[
\frac{\|m_{p,q}(x_m)-m_{p,q}(y_m)\|_q}
{\|x_m-y_m\|_p^\alpha}
=
\frac{\|m_{p,q}(x_m)-m_{p,q}(y_m)\|_q}
{\|x_m-y_m\|_p^{p/q}}
\cdot
\frac1{\|x_m-y_m\|_p^{\beta p/q}}
\longrightarrow\infty.
\]
Hence $\frac pq$ is the optimal Hölder exponent.

\smallskip

\noindent
{\it Optimality of the Hölder constant.}

Let
\[
x=
(2^{-1/p},0,\dots,0,2^{-1/p}),
\qquad
y=
(-2^{-1/p},0,\dots,0,2^{-1/p}).
\]
Then
\[
\|m_{p,q}(x)-m_{p,q}(y)\|_q
=
2^{1-\frac pq}
\|x-y\|_p^{p/q},
\]
showing that
$2^{1-\frac pq}$ is the optimal Hölder constant.

\medskip

\noindent
{\bf Case 2.}
Suppose $1\le q<p<\infty$.

\smallskip

\noindent
{\it Optimality of the exponent.}

The argument in the proof of Lemma~\ref{Q1} remains valid in the present finite-dimensional setting and hence shows that $m_{p,q}$ cannot be $\alpha$-H\"older for any $\alpha>1$.
Hence the Lipschitz exponent is optimal.

\smallskip

\noindent
{\it Optimality of the Lipschitz constant.}

For sufficiently large $m$, let
\[
x_m=
\Big(
(1-\tfrac1m)^{1/p}2^{-1/p},
0,\dots,0,
(1+\tfrac1m)^{1/p}2^{-1/p}
\Big),
\]
\[
y_m=
\Big(
(1+\tfrac1m)^{1/p}2^{-1/p},
0,\dots,0,
(1-\tfrac1m)^{1/p}2^{-1/p}
\Big).
\]
A straightforward Taylor expansion yields
\[
\lim_{m\to\infty}
\frac{\|m_{p,q}(x_m)-m_{p,q}(y_m)\|_q}
{\|x_m-y_m\|_p}
=
\frac pq.
\]
Therefore $\frac pq$ is the optimal Lipschitz constant.

Finally,
$m_{p,q}^{-1}=m_{q,p}$.
Hence the corresponding lower estimates for
$\rho_{m_{p,q}}$
follow immediately, and all the estimates in
Theorem~\ref{th1.1} are sharp.
\hfill$\square$

\begin{thm}\label{th1.3}
 Let  $(\Omega,\mu)$ be a non-atomic probability space and let $1\leq p,q<\infty.$ Let $M_{p,q}: S_{L_p(\Omega,\mu,\mathbb{C})}\rightarrow S_{L_q(\Omega,\mu,\mathbb{C})}$ be the Mazur map defined by
\begin{equation*}
  M_{p,q}(f)=|f|^{\frac{p}{q}}{\rm sign}(f),~f\in S_{L_p(\Omega,\mu,\mathbb{C})}.
\end{equation*}
Then  those H\"{o}lder estimates appearing in Mazur's sphere theorem (Theorem \ref{th1.1}) are all sharp.

\end{thm}
\pf For the given probability space $(\Omega,\mu)$, since for each $n$ we may represent $\ell_p^n$ as the subspace of $L_p(\mu)$ consisting of functions that are constant on some fixed $n$ disjoint sets (and similarly for $\ell_q^n$)
, $M_{p,q}: S_{L_p(\Omega,\mu,\mathbb{C})}\rightarrow S_{L_q(\Omega,\mu,\mathbb{C})}$ restricts to the natural Mazur map $m_{p,q}$ from $S_{\ell_p^n}$ onto $S_{\ell_q^n}$. Since $m_{p,q}$ admits the sharp (H\"{o}lder) estimates  between $S_{\ell_p^n}$ and $S_{\ell_q^n}$,   $M_{p,q}$ cannot have better H\"{o}lder estimates  as a map from $S_{L_p(\Omega,\mu,\mathbb{C})}$ onto $S_{L_q(\Omega,\mu,\mathbb{C})}$ than the corresponding map from $S_{\ell^n_p}$ onto $S_{\ell^n_q}$. Thus  the H\"{o}lder estimates of  $M_{p,q}$   between $S_{L_p(\Omega,\mu,\mathbb{C})}$ and $S_{L_q(\Omega,\mu,\mathbb{C})}$ appearing in Theorem \ref{th1.1} are also  sharp.
\hfill$\square$
\subsection{Optimal H\"older exponents for sphere embeddings}

Mazur's sphere theorem provides explicit H\"older embeddings between
classical \(L_p\) spheres.  It is therefore natural to ask whether these
H\"older exponents are optimal among all uniform embeddings of spheres into
classical \(L_p\) spaces.  The following quantity formalizes this question.

\begin{defn}
Let \(\mathcal M\) and \(\mathcal N\) be metric spaces. Define
\(\omega_{\mathcal N}(\mathcal M)\) to be the infimum over all
\(\alpha>0\) with the following property: every uniform embedding
\(f:\mathcal M\to\mathcal N\) satisfies
\[
\omega_f(t)\ge C_f t^\alpha,\qquad 0\le t\le 2,
\]
for some constant \(C_f>0\).
\end{defn}

Although rather elementary, the above estimate serves as the starting
point of our sphere estimates  theory. Stronger estimates  can be obtained by
combining it with geometric properties of classical Banach spaces.

\begin{thm}\label{Sharp}
Let \(1\le p<q<\infty\). Then
\[
\omega_{\ell_q}(S_{\ell_p})=\frac pq,
\qquad
\omega_{\ell_p}(S_{\ell_q})=1,
\]
and
\[
\omega_{L_q}(S_{L_p})=
\begin{cases}
\dfrac pq, & 1\le p<q\le2,\\[10pt]
\dfrac p2, & 1\le p<2\le q,\\[10pt]
1, & 2\le p<q,
\end{cases}
\qquad
\omega_{L_p}(S_{L_q})=1.
\]
Moreover, in each case the above value is attained.
\end{thm}\pf
We first compute \(\omega_{\ell_q}(S_{\ell_p})\).  Since \(\ell_q\) is
reflexive and \(q\)-AUS, Theorems~\ref{KR2008} and~\ref{HFC} imply that
every uniform embedding $f: S_{\ell_p} \to \ell_q$
 satisfies
\(\omega_f(t)\gtrsim t^{p/q}\). Hence
\(\omega_{\ell_q}(S_{\ell_p})\le p/q\).
Conversely, the Mazur map
\(m_{p,q}:S_{\ell_p}\to S_{\ell_q}\subset\ell_q\)
satisfies \(\omega_{m_{p,q}}(t)\lesssim t^{p/q}\), so
\(\omega_{\ell_q}(S_{\ell_p})\ge p/q\). Therefore
\[
\omega_{\ell_q}(S_{\ell_p})=\frac pq.
\]

Next we consider $S_{\ell_q}\stackrel{\rm u}{\hookrightarrow}\ell_p
\quad\text{and}\quad
S_{L_q}\stackrel{\rm u}{\hookrightarrow}L_p.$
Let \(f:S_{\ell_q}\to\ell_p\) and \(g:S_{L_q}\to L_p\) be uniform embeddings.
By Lemma~\ref{Q1}, both satisfy \(\omega_f(t),\omega_g(t)\gtrsim t\).
Hence
\(\omega_{\ell_p}(S_{\ell_q})\le1\) and
\(\omega_{L_p}(S_{L_q})\le1\).
Conversely, the Mazur maps
\(m_{q,p}\) and \(M_{q,p}\) are Lipschitz, so
\(\omega_{\ell_p}(S_{\ell_q})\ge1\) and
\(\omega_{L_p}(S_{L_q})\ge1\). Thus
\[
\omega_{\ell_p}(S_{\ell_q})
=
\omega_{L_p}(S_{L_q})
=1.
\]

Finally, we consider $S_{L_p}\stackrel{\rm u}{\hookrightarrow}L_q.$
Let \(f:S_{L_p}\to L_q\) be a uniform embedding.

If \(2\le p<q\), Lemma~\ref{Q1} yields
\(\omega_f(t)\gtrsim t\), hence
\(\omega_{L_q}(S_{L_p})\le1\).

Assume now \(1\le p<q\le2\) or \(1\le p<2\le q\).
Since \(L_p\) contains an isometric copy of \(\ell_p\),
the restriction \(f|_{S_{\ell_p}}\) is a uniform embedding into \(L_q\).
Using that \(L_q\) is \(q\)-AUS for \(1<q<2\) and \(2\)-AUS for \(q\ge2\)
(see \cite{Milman1971,JLPS2002}),
Theorems~\ref{KR2008} and~\ref{HFC} give
\[
\omega_{L_q}(S_{L_p})\le
\begin{cases}
\dfrac pq,&1\le p<q\le2,\\[10pt]
\dfrac p2,&1\le p<2\le q.
\end{cases}
\]

For the lower bounds, if \(1\le p<q\le2\), the Mazur map
\(M_{p,q}:S_{L_p}\to S_{L_q}\) yields
\(\omega_{L_q}(S_{L_p})\ge p/q\).
If \(1\le p<2\le q\), let \(F=T\circ M_{p,2}\), where
\(T:L_2\to L_q\) is an isometric embedding.
Then \(\omega_F(t)\lesssim t^{p/2}\), hence
\(\omega_{L_q}(S_{L_p})\ge p/2\).
If \(2\le p<q\), the same construction gives
\(\omega_F(t)\lesssim t\), so
\(\omega_{L_q}(S_{L_p})\ge1\).
Combining these estimates yields
\[
\omega_{L_q}(S_{L_p})=
\begin{cases}
\dfrac pq,&1\le p<q\le2,\\[10pt]
\dfrac p2,&1\le p<2\le q,\\[10pt]
1,&2\le p<q.
\end{cases}
\]

The proof shows that all optimal exponents are attained by the Mazur map
or by its composition with an isometric embedding \(L_2\hookrightarrow L_q\).
\hfill$\square$

Theorem~\ref{Sharp} completely determines the optimal H\"older exponents
for uniform embeddings of classical \(L_p\)-spheres into \(L_q\)-spaces.
Moreover, these optimal exponents are always attained by the Mazur map,
or by its composition with an isometric embedding of \(L_2\) into
\(L_q\). Thus, despite its remarkably simple and explicit form, the
Mazur map provides an optimal model for uniform embeddings between
classical \(L_p\)-spheres.

This phenomenon parallels the nonlinear embedding theory of the ambient
\(L_p\)-spaces. Indeed, Mendel and Naor~\cite{MN2004} proved that for every
\(1\le q<p<\infty\), there exists an explicit uniform embedding
\[
f:L_q\longrightarrow L_p
\]
such that
\[
\|f(x)-f(y)\|_p
=
\|x-y\|_q^{\,q/p},
\qquad x,y\in L_q,
\]
and then Naor demonstrated that the exponent \(q/p\) is optimal for $2<q<p<\infty$\cite{Naor2016}.

\subsection{Sphere Version of Naor's bi-$\theta$-Hölder  embedding Theorem }

\quad
A remarkable result of Mendel and Naor \cite{MN2004} shows that  for every $1\le q<p<\infty$ there exists an explicit mapping $f:L_q\longrightarrow L_p$
satisfying
\[
\|f(x)-f(y)\|_{L_p}
=\|x-y\|_{L_q}^{q/p},
\qquad x,y\in L_q.
\]
Consequently $L_q$ admits a bi-$\frac qp$-Hölder embedding into $L_p$.
A fundamental question asks whether the exponent
\(\frac{q}{p}\) is optimal. Metric type and cotype yield sharp restrictions
on such exponents when either $1\le q\le2$ and $p\ge q$,
or $2\le q<\infty$ and $p\le q$. For $2<q<p<\infty$, Naor and Schechtman
(Conjecture~1.8 in \cite{NS2016}) conjectured that this exponent is
optimal: any bi-$\theta$-H\"older embedding of $L_q$
into $L_p$ must satisfy $\theta\le q/p$.
This conjecture was settled affirmatively by Naor \cite{Naor2016}, who proved that, in this range, the maximal admissible
exponent $\theta\in(0,1]$ equals \(\frac{q}{p}\).
A key ingredient in the proof is the sharp metric \(X_p\) inequality,
namely that \(L_p\) is an \(X_p\) metric space with the sharp scaling
parameter.

Restricting the Mendel--Naor embedding to the unit sphere yields a
bi-\(\frac{q}{p}\)-H\"older embedding of the unit sphere $S_{L_q}$ into $L_p$.
It is therefore natural to ask whether the exponent
\(\frac{q}{p}\) remains optimal in the sphere setting.
The following theorem gives a positive answer.
Our proof combines Naor’s sharp metric \(X_p\) inequality \cite{Naor2016} with suitable maps from discrete tori into \(S_{L_q}\) to obtain the required estimates on the unit sphere.

\begin{thm}[Naor \cite{Naor2016}]
Let  \(2<q<p<\infty\). Suppose that there exist
\(\theta\in(0,1]\), \(L\ge1\), and $f:L_q\to   L_p$ such that
\[
\|x-y\|_{L_q}^{\theta}
\le
\|f(x)-f(y)\|_{L_p}
\le
L\|x-y\|_{L_q}^{\theta},
\qquad x,y\in L_q.
\]
Then $\theta\le \frac{q}{p}.$

\end{thm}

The following theorem shows that the same sharp exponent remains valid
for unit spheres. 
\begin{thm}[Sphere version of Naor's bi-\(\theta\)-H\"older embedding theorem]\label{sphere-holder-critical} 
Let \(2<q<p<\infty\). Suppose that there exist \(\theta\in(0,1]\),
\(L\ge1\), and  $f: S_{L_q}\to L_p$ such that
\[
\|x-y\|_{L_q}^{\theta}
\le
\|f(x)-f(y)\|_{L_p}
\le
L\|x-y\|_{L_q}^{\theta},
\qquad x,y\in S_{L_q}.
\]
Then $\theta\le \frac{q}{p}.$

\end{thm}
\pf
Assume, toward a contradiction, that $\theta>\frac qp.$
Then
\[
\rho_f(t)\ge t^\theta,
\qquad
\omega_f(t)\le Lt^\theta .
\]

By Theorem~\ref{sphere-modulus-ineq}, for every \(n\in\mathbb N\),
\(k\in[n]\), and \(m\ge \sqrt{n/k}\),
\[
\rho_f\!\left(
2\left(\frac{k}{n}\right)^{1/q}
\right)
\le
C\, m
\left(
k\,\omega_f\!\left(\frac{\pi}{4m}n^{-1/q}\right)^p
+
\left(\frac{k}{n}\right)^{p/2}
\omega_f\!\left(\frac{\pi}{4m}\right)^p
\right)^{1/p}.
\]

Using the Hölder estimates for \(f\), we obtain
\begin{align*}
\left(2\left(\frac{k}{n}\right)^{1/q}\right)^\theta
&\le
C\, m
\left(
k\,L^p\left(\frac{\pi}{4m}n^{-1/q}\right)^{\theta p}
+
\left(\frac{k}{n}\right)^{p/2}
L^p\left(\frac{\pi}{4m}\right)^{\theta p}
\right)^{1/p} \\
&\le
C L\, m^{1-\theta}
\left(
k^{1/p}n^{-\theta/q}
+
\left(\frac{k}{n}\right)^{1/2}
\right).
\end{align*}

Hence
\[
\left(\frac{k}{n}\right)^{\theta/q}
\le
C L\, m^{1-\theta}
\left(
k^{1/p}n^{-\theta/q}
+
\left(\frac{k}{n}\right)^{1/2}
\right).
\]

Choose \(a\in(0,1)\) sufficiently close to \(1\) such that
\[
(1-a)\frac{1-\theta}{2}
+
a\left(\frac1p-\frac{\theta}{q}\right)
<0,
\]
which is possible since \(\frac1p-\frac{\theta}{q}<0\).
Let
\[
k=\lceil n^a\rceil,
\qquad
m=\left\lceil \sqrt{\frac nk}\right\rceil .
\]
Here and below, \(\lceil t\rceil\) denotes the smallest integer greater than
or equal to \(t\). 
Then
\[
k\asymp n^a,
\qquad
m\asymp n^{(1-a)/2}.
\]

Substituting into the previous inequality yields
\[
n^{-(1-a)\theta/q}
\le
CL
\left(
n^{(1-a)(1-\theta)/2+a/p-\theta/q}
+
n^{-(1-a)\theta/2}
\right).
\]

Dividing by \(n^{-(1-a)\theta/q}\), we get
\[
1
\le
CL
\left(
n^{(1-a)(1-\theta)/2+a/p-a\theta/q}
+
n^{-(1-a)\theta(1/2-1/q)}
\right).
\]

By the choice of \(a\), both exponents are negative. Letting \(n\to\infty\)
yields a contradiction. Therefore $\theta\le \frac qp .$
\hfill$\square$

The above result concerns Hölder embeddings. We now turn to a stronger
quantitative category, namely almost Lipschitz embeddings, and show that
our estimates also yield rigidity phenomena beyond Hölder geometry.

\section{Applications to the Almost Lipschitz Geometry of $L_p$ Spheres}
\subsection{Almost Lipschitz Maps and Embeddings}

\quad

Mazur’s sphere theorem shows that the unit spheres of separable
infinite-dimensional $L_p$ spaces are uniformly homeomorphic and even
Hölder homeomorphic.

Thus, neither uniform nor Hölder structures are sufficient to
distinguish the geometry of $L_p$ spheres.

On the other hand, Lipschitz embeddings are strong enough to capture
quantitative geometric features, but are often too rigid to allow
natural nonlinear constructions.

This leads to a natural question: is there a notion of embedding whose
modulus of continuity lies strictly between Hölder and Lipschitz
behaviour? More precisely, can one identify a class of embeddings that
is stronger than Hölder control, yet still weaker than Lipschitz
regularity, and is sufficiently flexible while retaining asymptotic
Lipschitz behaviour at small scales?

This motivates the notion of almost Lipschitz maps introduced below, which has also appeared recently in \cite{xiang2026}.

\begin{defn}[Almost Lipschitz maps and embeddings]\label{ALIP}
Let $(\mathcal M,d_{\mathcal M})$ and
$(\mathcal N,d_{\mathcal N})$ be metric spaces, and let
$f:\mathcal M\to\mathcal N$ be a mapping.

\begin{enumerate}
\item[(i)]
Let $\alpha>0$. We say that $f$ is
\emph{$\alpha$-almost Lipschitz} if there exist constants
$A>0$ and $t_0>0$ such that
\[
\omega_f(t)\le At|\log t|^\alpha,
\qquad
0<t\le t_0.
\]

\item[(ii)]
Let $\alpha,\beta>0$. An injective mapping
$f:\mathcal M\to\mathcal N$
is called an
\emph{$(\alpha,\beta)$-almost Lipschitz embedding}
if $f$ is $\alpha$-almost Lipschitz and
$f^{-1}:f(\mathcal M)\to\mathcal M$
is $\beta$-almost Lipschitz.

\item[(iii)]
We say that
$\mathcal M$ is
\emph{almost Lipschitz embeddable}
into $\mathcal N$, denoted by
\[
\mathcal M
\stackrel{\rm AL}{\hookrightarrow}
\mathcal N,
\]
if there exists an
$(\alpha,\beta)$-almost Lipschitz embedding
from $\mathcal M$ into $\mathcal N$
for some
$\alpha,\beta>0$.

\item[(iv)]
A bijection
$f:\mathcal M\to\mathcal N$
is called an
\emph{$(\alpha,\beta)$-almost Lipschitz homeomorphism}
if both $f$ and $f^{-1}$ are almost Lipschitz with exponents
$\alpha$ and $\beta$, respectively.

The metric spaces
$\mathcal M$
and
$\mathcal N$
are said to be
\emph{$(\alpha,\beta)$-almost Lipschitz homeomorphic}
if there exists an
$(\alpha,\beta)$-almost Lipschitz homeomorphism between them.
\end{enumerate}
\end{defn}

The above definition should be viewed as a quantitative refinement of
uniform geometry. It naturally interpolates between Lipschitz and
uniform embeddings:
\[
\text{\rm Lipschitz}
\Longrightarrow
\text{\rm almost Lipschitz}
\Longrightarrow
\text{\rm uniform},
\]
and both implications are strict in general.


\begin{rem}\label{rem:rho}
For an injective map \(f\), estimates on \(\omega_{f^{-1}}\) may be
equivalently expressed in terms of the compression modulus
\(\rho_f\). In particular, for every \(\beta>0\),
\[
\omega_{f^{-1}}(t)\lesssim
t|\log t|^\beta
\quad\Longleftrightarrow\quad
\rho_f(t)\gtrsim
\frac{t}{|\log t|^\beta},
\]
for all sufficiently small \(t>0\). Accordingly, an
\((\alpha,\beta)\)-almost Lipschitz embedding can be formulated using
the pair \((\omega_f,\rho_f)\), and we shall adopt this formulation
throughout the paper.
\end{rem}

\begin{rem}

\quad

\begin{enumerate}

\item[(i)]
Every Lipschitz mapping is $\alpha$-almost Lipschitz for every
$\alpha>0$, while every almost Lipschitz mapping is uniformly
continuous. Neither converse holds in general.

\item[(ii)]
The idea underlying almost Lipschitz embeddings can already be found in
Benyamini--Lindenstrauss
\cite[Proposition~7.18]{BL2000}, who proved that every compact metric
space $\mathcal M$ admits a uniform embedding $f:\mathcal M\longrightarrow
\Big(\sum_{n=1}^{\infty}\ell_\infty^n\Big)_{\ell_2}$
such that
\[
\omega_f(t)\lesssim t|\log t|^{1/2},
\qquad
\rho_f(t)\gtrsim t.
\]
This provides a prototypical example of an almost Lipschitz embedding
in the sense of the present paper.

\item[(iii)]
On the other hand, \cite[Proposition~4.1(iii)]{BL2015} shows that there exist compact subsets of $c_0$
which do not admit Lipschitz embeddings into $\Big(\sum_{n=1}^{\infty}\ell_{\infty}^{n}\Big)_{\ell_{2}}.$
We now consider the compact set
\[
K=\overline{\mathrm{co}}\Big(\Big\{\frac{e_n}{n}\Big\}_{n=1}^{\infty}\cup\{0\}\Big)\subset c_0.
\]
We claim that $K$ does not admit a Lipschitz embedding into
$\big(\sum_{n=1}^{\infty}\ell_{\infty}^{n}\big)_{\ell_{2}}$.

Indeed, if such an embedding existed, then by
\cite[Theorem~3.6]{CZ2009}
the space $c_0$ would embed linearly into
$\big(\sum_{n=1}^{\infty}\ell_{\infty}^{n}\big)_{\ell_{2}}$,
which is impossible.

\item[(iv)]
The notion introduced above should not be confused with the
\emph{almost Lipschitz embeddability} developed by
Baudier and Lancien \cite{BL2015,B2022}. Their notion is formulated in
terms of a family of embeddings whose compression behaviour can be made
arbitrarily close to linear at different scales, whereas ours is defined
directly through the moduli of continuity of a single embedding.

Consequently, the two notions reflect different quantitative structures
and are not directly comparable.

Our formulation is particularly suited to the study of individual
uniform embeddings and sphere mappings, since it allows one to work
directly with the pair
\[
(\omega_f,\omega_{f^{-1}}),
\quad\text{or equivalently}\quad
(\omega_f,\rho_f).
\]

\item[(v)]
The terminology ``almost Lipschitz'' is motivated by two complementary
facts. First, it already appears implicitly in Benyamini--Lindenstrauss
\cite[Proposition~7.18]{BL2000}, where embeddings with a logarithmic perturbation of Lipschitz regularity are constructed and described as
``almost Lipschitz''. Second, from a quantitative viewpoint, an
$\alpha$-almost Lipschitz map is $\beta$-Hölder for every
$0<\beta<1$ and for sufficiently small distances, with estimates improving as $\beta\nearrow 1$. Thus such
mappings approach the Lipschitz scale arbitrarily closely while still
allowing a mild logarithmic loss.

These two viewpoints together justify the terminology: the class of
almost Lipschitz mappings lies quantitatively between Hölder and
Lipschitz regularity.

\end{enumerate}
\end{rem}

\subsection{Sphere version of the Enflo-Lindenstrauss theorem} 

The following theorem is the classical result of Lindenstrauss \cite{Lindenstrauss1964} and Enflo \cite{Enflo1969}.
Our next theorem establishes its sphere analogue in the almost Lipschitz
category.
\begin{thm}[Enflo--Lindenstrauss]
Let $1\le p,q<\infty$, and let
$L_p(\mu_1)$ and $L_q(\mu_2)$
be separable Banach spaces.
If $L_p(\mu_1)$ and $L_q(\mu_2)$ are uniformly homeomorphic, then either both spaces are finite-dimensional and have the same dimension, or $p=q$.
\end{thm}

\begin{thm}[Sphere version of the Enflo--Lindenstrauss theorem]\label{sphere-Enflo-Lindenstrauss}
Let $1\le p,q<\infty$, and let
$L_p(\mu_1)$ and $L_q(\mu_2)$
be separable Banach spaces.
If the unit spheres
$S_{L_p(\mu_1)}$ and $S_{L_q(\mu_2)}$
are almost Lipschitz homeomorphic, then either both spaces are finite-dimensional and have the same dimension, or $p=q$.
\end{thm}
\pf
If one of the two spaces is finite-dimensional, then so is the other, since
the unit sphere of a Banach space is compact if and only if the space itself
is finite-dimensional. As the two spheres are homeomorphic, the invariance
of domain implies that the two spaces have the same dimension.

Assume therefore that both spaces are separable and infinite-dimensional,
and suppose, toward a contradiction, that $p<q$. Let
\[
f:S_{L_p(\mu_1)}\longrightarrow S_{L_q(\mu_2)}
\]
be an almost Lipschitz homeomorphism. By the
classification of separable infinite-dimensional $L_r$-spaces, each such
space is isomorphic to either $\ell_r$ or $L_r[0,1]$; see
\cite[p.~137]{AK2016}. Thus the situation reduces to the four pairs
\[
(\ell_p,\ell_q),\qquad
(\ell_p,L_q),\qquad
(L_p,\ell_q),\qquad
(L_p,L_q).
\]
Assume first that $1\le p<2$. Combining this with Theorem \ref{KR2008}, Theorem~\ref{HFC}, and the fact that $L_r$ is $r$-AUS for
$1<r\le2$ and $2$-AUS for $2\le r<\infty$
\cite{Milman1971,JLPS2002}, we obtain
\[
\omega_f(t)\gtrsim t^\theta
\]
for some $\theta=\theta(p,q)\in(0,1)$. On the other hand, the almost
Lipschitz assumption yields
\[
\omega_f(t)\lesssim t|\log t|^\alpha
\]
for some $\alpha>0$. This is impossible, since
$t|\log t|^\alpha=o(t^\theta)$ as $t\to0^+$ for every $\alpha>0$ and
$0<\theta<1$. This rules out the case $1\le p<2$.

Assume now that $2\le p<q<\infty$.
After the above reduction, write
\[
X\in\{\ell_p,L_p\},
\qquad
Y\in\{\ell_q,L_q\},
\]
and regard $f$ as an almost Lipschitz homeomorphism
from $S_X$ onto $S_Y$.
Let $J:\ell_q\longrightarrow Y$ be a linear isometric
embedding, and define
\[
g=f^{-1}\circ J\big|_{S_{\ell_q}}:
S_{\ell_q}\longrightarrow X.
\]
Then $g$ is an almost Lipschitz embedding.
Thus, there exist
$\alpha,\beta>0$ such that, for all sufficiently small $t>0$,
\[
\omega_g(t)\lesssim t|\log t|^\alpha,
\qquad
\rho_g(t)\gtrsim \frac{t}{|\log t|^\beta}.
\]
Since $X$ has cotype $p$ and nontrivial type, it follows from
Theorem~\ref{Cotype1}  that
\[
\rho_g(t)\lesssim\omega_g(t^{q/p}).
\]
Therefore,
\[
\frac{t}{|\log t|^\beta}
\lesssim \rho_g(t)
\lesssim \omega_g(t^{q/p})
\lesssim t^{q/p}|\log t|^\alpha.
\]
This implies that
\[
1\lesssim t^{q/p-1}|\log t|^{\alpha+\beta}
\longrightarrow0
\qquad (t\to0^+),
\]
a contradiction, since $q/p>1$.
Thus the case $2\le p<q$ is also impossible,
and therefore $p<q$ cannot occur.

Applying the same argument to $f^{-1}$ rules out the case $p>q$. Therefore
$p=q$.
\hfill$\square$

\begin{rem}

\quad

\begin{enumerate}

\item[(i)]
Mazur's sphere theorem asserts that the unit spheres of separable
infinite-dimensional \(L_p(\mu)\) spaces are uniformly homeomorphic.
Moreover, since the Mazur map and its inverse are Hölder continuous,
these spheres are in fact Hölder homeomorphic.
Theorem~\ref{sphere-Enflo-Lindenstrauss} shows that this phenomenon no
longer persists in the almost Lipschitz category: two separable
\(L_p(\mu)\)-spheres are almost Lipschitz homeomorphic only if they
arise from the same exponent \(p\) (apart from the trivial
finite-dimensional case). Hence almost Lipschitz homeomorphisms provide
a strictly finer classification than both uniform and Hölder
homeomorphisms for \(L_p (\mu)\) spheres.

\item[(ii)]
The above theorem illustrates the central theme of this paper: the
quantitative behaviour of the moduli of continuity provides a finer
geometric invariant than uniform or Hölder equivalence. In particular,
it captures geometric features of \(L_p\)-spheres that are invisible to
the classical uniform and Hölder categories.

\end{enumerate}
\end{rem}





\subsection{Sphere Version of the Banach--Kadec--Paley Theorem}

\quad

The classical Banach--Kadec--Paley theorem (\cite{Banach1932,Kadec1958,Paley1936}) characterizes when one classical \(L_q\)-space embeds linearly into another \(L_p\)-space. Using the notation introduced earlier, we write
$X\hookrightarrow Y$ and $X\stackrel{\rm AL}{\hookrightarrow}Y$
to denote a linear embedding and an almost Lipschitz embedding of $X$ into $Y$,
respectively.

\begin{thm}[Banach--Kadec--Paley theorem]
Let \(1\le p\neq q<\infty\). Then 
\[
L_q\hookrightarrow L_p \quad\Longleftrightarrow\quad q=2 \quad {\rm or} \quad 1\le p< q\le2.
\]
\end{thm}

The following theorem is the main result of this paper. It shows that,
for \(L_p\)-spaces, the almost Lipschitz geometry of unit spheres
completely determines the  linear structure of the space. 

\begin{thm}[Sphere version of the Banach--Kadec--Paley theorem]
\label{thm:spherical-KPP}

Let \(1\le p,q<\infty\). Then \[
S_{L_q}\stackrel{\rm AL}{\hookrightarrow}L_p
\quad\Longleftrightarrow\quad
L_q\hookrightarrow L_p.
\]
\end{thm}
\pf
Suppose $S_{L_q}\stackrel{\rm AL}{\hookrightarrow}L_p,$
and let
\(f:S_{L_q}\to L_p\)
be an almost Lipschitz embedding. Then there exist
\(\alpha,\beta>0\) such that, for all sufficiently small \(t>0\),
\[
\omega_f(t)\lesssim
t\Bigl(\log\frac1t\Bigr)^\alpha,
\qquad
\rho_f(t)\gtrsim
\frac{t}{\bigl(\log\frac1t\bigr)^\beta}.
\]

 We distinguish the following
three cases according to the relative positions of \(p\) and \(q\).

\medskip
\noindent
\textbf{Case 1.} \(q<\min\{p,2\}\).

Since \(L_q\) contains an isometric copy of \(\ell_q\), the restriction
of \(f\) to \(S_{\ell_q}\) is still an almost Lipschitz embedding,
which we continue to denote by \(f\).
By Theorem~\ref{HFC},
\[
\omega_f(\varepsilon)\gtrsim
\begin{cases}
\varepsilon^{q/p}, & 1\le q<p\le2,\\[2mm]
\varepsilon^{q/2}, & 1\le q<2\le p.
\end{cases}
\]
Since both exponents are strictly smaller than \(1\), this contradicts
\[
\omega_f(\varepsilon)
\lesssim
\varepsilon
\Bigl(\log\frac1\varepsilon\Bigr)^\alpha .
\]

\medskip
\noindent
\textbf{Case 2.} \(2<q<p\).

 By Theorem~\ref{sphere-modulus-ineq}, there exists a
constant \(C_p>0\), depending only on \(p\), such that for every
\(n\in\mathbb N\), every \(k\in[n]\), and every
\(m\ge\sqrt{n/k}\),
\[
\rho_f\!\left(2\Bigl(\frac{k}{n}\Bigr)^{1/q}\right)
\le
C_pmk^{1/p}
\omega_f\!\left(\frac{\pi}{4m}n^{-1/q}\right)
+
C_pm\Bigl(\frac{k}{n}\Bigr)^{1/2}
\omega_f\!\left(\frac{\pi}{4m}\right).
\]

 Choose $m, k \in \mathbb{N}$ as follows.
\[
m \stackrel{\mathrm{def}}{=} \left\lfloor n^{\frac{ p - q}{q(p-2)}} \right\rfloor \quad \text{and} \quad k \stackrel{\mathrm{def}}{=} \left\lceil \dfrac{ n}{m^2} \right\rceil.
\]
Then we have $m\leq n^{\frac{ p - q}{q(p-2)}} \leq 2m$ and 
$1\leq \frac{ n}{m^2}\leq k\leq \frac{2n}{m^2}.$
Clearly, $k\in \{1,2,\cdots,n\}$ and $ m \geqslant \sqrt{n/k}$. Thus, we have

\[
\rho_{f}\left(2\left(\frac{k}{n}\right)^{\frac{1}{q}}\right)
\leq C_{p}\left(2^{\frac{1}{p}}n^{\frac{1}{q}}\omega_{f}\left(\frac{n^{-\frac{1}{q}}\pi}{4m}\right)+\sqrt{2}\omega_{f}\left(\frac{\pi}{4m}\right)\right)\]

Now, for large enough $m$ we have
$\frac{1}{q}\log n \leq \frac{2(p-2)}{p-q}\log m$
and hence 
 \begin{eqnarray*}
\rho_{f}\left(2 m^{-\frac{2}{q}}\right)&\leq& \rho_{f}\left(2\left(\frac{k}{n}\right)^{\frac{1}{q}}\right)\\
&\leq& C_{p}2^{1+\frac{1}{p}}n^{\frac{1}{q}}\omega_{f}\left(\frac{n^{-\frac{1}{q}}}{m}\right)+ 2\sqrt{2}C_{p}\omega_{f}\left(\frac{1}{m}\right)\\
&\lesssim&\frac{1}{m}(\log m)^{\alpha}
 \end{eqnarray*}
 Moreover, we obtain
 \[\rho_{f}\left(2 m^{-\frac{2}{q}}\right)\gtrsim \dfrac{m^{-\frac{2}{q}}}{(\log m)^{\beta}}.\]

 This implies that
 \[1\lesssim m^{\frac{2}{q}-1}(\log m)^{\alpha+\beta}.\]
 Since $q>2$, letting $n\to\infty$ (hence $m\to\infty$), this leads to a contradiction.

\medskip
\noindent
\textbf{Case 3.} \(q>\max\{p,2\}\).

Assume first that \(p>1\), and let \(r=\max\{p,2\}\).
As in Case 1, the restriction of \(f\) to  \(S_{\ell_q}\) is still denoted by \(f\). 
Since \(L_p\) has nontrivial type and cotype \(r\), it follows from 
Theorem~\ref{Cotype1} that
\[
\rho_f(t)\lesssim\omega_f(t^{q/r}).
\]
Thus, 
\[
\frac{t}{\bigl(\log(1/t)\bigr)^\beta}
\lesssim
t^{q/r}\left(\log\frac{1}{t}\right)^\alpha.
\]
Consequently,
\[
1
\lesssim
t^{q/r-1}\left(\log\frac{1}{t}\right)^{\alpha+\beta}
\longrightarrow 0
\qquad \text{as } t\to 0^+,
\]
because \(q/r>1\), which is a contradiction.

It remains to consider the endpoint \(p=1\). By
\cite[Remark~5.10]{MN2004}, there exists a uniform embedding
\[
\varphi:L_1\longrightarrow L_2
\]
satisfying
\[
\|\varphi(x)-\varphi(y)\|_2
=\|x-y\|_1^{1/2},
\qquad x,y\in L_1.
\]
Hence
\(
\rho_\varphi(s)\ge s^{1/2}
\)
and
\(
\omega_\varphi(s)\le s^{1/2}.
\)

Applying Theorem~\ref{Cotype1} to
\(g=\varphi\circ f:S_{\ell_q}\to L_2\)
yields
\(
\rho_g(t)\lesssim\omega_g(t^{q/2}),
\)
while
\[
\rho_g(t)\ge\rho_f(t)^{1/2},
\qquad
\omega_g(t)\le\omega_f(t)^{1/2}.
\]
Therefore
\[
\rho_f(t)\lesssim\omega_f(t^{q/2}).
\]
 Since \(q/2>1\), the same argument as in the preceding \(p>1\) case
yields a contradiction.

Thus
\[
L_q\hookrightarrow L_p.
\]

Conversely,
\[
L_q\hookrightarrow L_p
\Longrightarrow
S_{L_q}\stackrel{\rm AL}{\hookrightarrow}L_p
\]
is obvious.
\hfill\(\square\)

\subsection{Almost Lipschitz homeomorphisms between Banach spaces}

\quad

In this section, we aim to  construct  a pair of separable Banach spaces that are almost Lipschitz homeomorphic but not Lipschitz homeomorphic. 

The example is derived from classical construction methods (see a very recent excellent book of Baudier and Lancien \cite[Chapter 13]{BL2026}). Specifically, we first present the basic lifting method  of almost Lipschitz homeomorphisms, and then combine this method with Lipschitz-free spaces to generate the desired example.

We begin with the lifting method. Suppose  $Q: X\to{Y}$ is a quotient map.  A \emph{(global) uniformly continuous lifting} of $Q$ is a uniformly continuous map $\varphi:Y\to{X}$ such that $Q\circ\varphi=id_{Y}$.  In particular, we say that $Q$ has a \emph{(global) $\alpha$-almost Lipschitz lifting} if the uniformly continuous lifting $\varphi: Y\to{X}$ is
 an $\alpha$-almost Lipschitz map for some $\alpha>0$.
 \begin{proposition}\label{pr}
 Let $Q: X\to{Y}$ be a canonical quotient map. If $Q$ admits an $\alpha$-almost Lipschitz lifting $\varphi$ for some $\alpha>0$, then $X$ and $Y\oplus_{1} {\rm ker} Q$ are
 $(\alpha,\alpha)$-almost Lipschitz homeomorphic.
 \end{proposition}
 \pf 
Define
\[
\Phi:Y\oplus_1 \ker Q\longrightarrow X,\qquad
\Phi(y,z)=\varphi(y)+z .
\]
Then \(\Phi\) is bijective and
\[
\Phi^{-1}(x)=\big(Qx,\,x-\varphi(Qx)\big).
\]
Since \(\varphi\) is \(\alpha\)-almost Lipschitz, there exist \(K>0\) and
\(t_0>0\) such that
\[
\omega_\varphi(t)\le Kt|\log t|^\alpha,\qquad 0<t\le t_0 .
\]
We may assume that \(t_0\le e^{-\alpha-1}\), so that
\(t\mapsto t|\log t|^\alpha\) is increasing on \((0,t_0)\).

Let \((y,z),(y',z')\in Y\oplus_1\ker Q\) with
\(\|(y,z)-(y',z')\|\le t_0\). Then
\[
\begin{aligned}
\|\Phi(y,z)-\Phi(y',z')\|
&\le \|\varphi(y)-\varphi(y')\|+\|z-z'\| \\
&\le K\|y-y'\||\log\|y-y'\||^\alpha
      +\|z-z'\||\log\|z-z'\||^\alpha \\
&\le (K+1)\|(y,z)-(y',z')\|
        |\log \|(y,z)-(y',z')\||^\alpha .
\end{aligned}
\]
Thus \(\Phi\) is \(\alpha\)-almost Lipschitz.

Similarly, if \(x,x'\in X\) and \(\|x-x'\|\le t_0\), then
\[
\begin{aligned}
\|\Phi^{-1}(x)-\Phi^{-1}(x')\|
&=\|Qx-Qx'\|
  +\|(x-x')-(\varphi(Qx)-\varphi(Qx'))\| \\
&\le 2\|x-x'\|+\|\varphi(Qx)-\varphi(Qx')\| \\
&\le (K+2)\|x-x'\||\log\|x-x'\||^\alpha .
\end{aligned}
\]
Hence \(\Phi^{-1}\) is also \(\alpha\)-almost Lipschitz. Therefore
\(\Phi\) is an \((\alpha,\alpha)\)-almost Lipschitz homeomorphism between
\(Y\oplus_1\ker Q\) and \(X\).
 \hfill$\square$
\par

Our examples relate closely to Lipschitz free spaces and we refer to \cite[\text{ pp. }178\text{-}188]{Kalton2004} for an introduction to Lipschitz free spaces.

Recall a \emph{gauge} is  a function $\omega: [0,\infty)\to [0,\infty)$ which is a continuous
increasing subadditive function with $\omega(0)=0$ and $\omega(t) \geq  t$ for $0 \leq  t \leq 1$. We say $\omega$ is
normalized if $\omega(1)=1$ and nontrivial if $\lim\limits_{t\to 0^+}\frac{\omega(t)}{t}=\infty.$ Moreover, a gauge $\omega$ is  called strongly normalized if  $\omega(t)=t$ for all
$t\geq 1$.  Such gauges do not distort the metric at
large distances and are appropriate for the study of almost Lipschitz homeomorphisms between
Banach spaces.

Let $X$ be a Banach space and denote metric $d$ by $d(x,y)=\|x-y\|$.  For any gauge $\omega$, the space ${\rm Lip}_{\omega}(X)$ of
$\omega\circ{d}$-Lipschitz functions on $X$ which vanish at $0$ has a natural predual, denoted $\mathcal{F}_{\omega}(X)$, which is the Lipschitz-free space over $X$. Aliaga, Gartland, Petitjean and Proch\'{a}zka \cite[after Theorem 4.6]{AGPP2022} showed that if $\omega$  is a nontrivial gauge then $\mathcal{F}_{\omega}(X)$ has the Radon-Nikod\'ym property.

We recall that,  if $\omega$ is a strongly normalized gauge,  there is a natural quotient map (the barycentric map) $\beta_{\omega}: \mathcal{F}_{\omega}(X)\to X$. Moreover,  $\delta_{\omega}: X\to\mathcal{F}_{\omega}(X)$ is uniformly continuous lifting  of $\beta_{\omega}$ which satisfies
\[\|\delta_{\omega}(x)-\delta_{\omega}(y)\|=\omega(\|x-y\|), \,\,x,y\in{X}.\]

By applying Proposition \ref{pr}, the following example of a strongly normalized gauge is at the heart of the construction of examples of almost Lipschitz homeomorphic Banach spaces.

\begin{proposition}\label{PR1}
Given $0<\alpha\leq 2$, let $\omega: [0,\infty)\to[0,\infty)$ be a function defined by
\[\omega(t)=\begin{cases}0,&t=0\\
t(|\log t|^\alpha+1),&0<t<e^{-1}\\
\frac{e-2}{e-1}t+\frac{1}{e-1},&e^{-1}\leq t< 1\\
t,&t\geq 1.\\

\end{cases}\]
Then $\omega:[0,\infty)\to [0,\infty)$ is a  nontrivial and strongly normalized gauge. That is, $\omega$ admits the following properties:

1) $\omega(0)=0,\omega(1)=1$.

2) $\omega(t)\geq t$ for $t\in [0,1];$ $\omega(t)=t$ for $t\in [1,\infty).$

3) $\lim\limits_{t\to 0^+}\frac{\omega(t)}{t}=\infty.$

4) $\omega:[0,\infty)\to [0,\infty)$ is continuous.

5) $\omega:[0,\infty)\to [0,\infty)$ is increasing.

6) $\omega:[0,\infty)\to [0,\infty)$ is subadditive.

\end{proposition}
\pf

Part 1) is obvious.

Part 2): For $0<t<e^{-1}$ we clearly have $\omega(t)=t(|\log t|^\alpha+1)>t.$ For $e^{-1}\leq t< 1$, since $\frac{e-2}{e-1}t+\frac{1}{e-1}-t=\frac{1-t}{e-1}>0$ we have $\omega(t)>t.$
Thus $\omega(t)\geq t$ for $t\in [0,\infty).$

Part 3) holds true since $\lim\limits_{t\to 0^+}\frac{\omega(t)}{t}=\lim\limits_{t\to 0^+}\frac{t(|\log t|^\alpha+1)}{t}=\infty.$

It remains to prove arguments 4),5) and 6).

First for 4): Clearly, $\omega$ is continuous on  $[0,\infty).$

Part 5) $\omega$ is increasing on $[0,\infty):$ Clearly $\omega$ is strictly increasing on $(e^{-1},1)$;  and $\omega$ is strictly increasing on $ (1,\infty)$.
Let $t\in (0,e^{-1})$. Then \begin{eqnarray*}
                                  \omega^\prime(t) &=& \Big[t \left((-\log t)^\alpha+1\right)\Big]^\prime \\
                                   &=& \Big[(-\log t)^\alpha+1\Big]+t\Big[\alpha(-\log t)^{\alpha-1}\cdot (-\frac{1}{t})\Big]\\
                                  &=& (-\log t)^\alpha-\alpha(-\log t)^{\alpha-1}+1. \\
                                  \end{eqnarray*}
Let $u=-\log t.$ Then $u>1$ for $0<t<e^{-1}$. We have \[\omega^\prime(t)=g(u)=u^\alpha-\alpha u^{\alpha-1}+1.\] Note that
                                   \[g^\prime(u)=\alpha u^{\alpha-1}-\alpha(\alpha-1)u^{\alpha-2}=\alpha u^{\alpha-2}(u-(\alpha-1)).\]Since $0<\alpha \leq 2$ and $u>1$ we have $g^\prime(u)>0.$ So $g$ is strictly increasing for $u>1$. Thus $g(u)>g(1)=2-\alpha\geq 0.$ Thus \[\omega^\prime(t)>0, ~0<t<e^{-1}.\]
                                  We have that $\omega$ is strictly increasing on $(0,e^{-1}).$

                                  Since $\omega$ is continuous on $[0,\infty)$ we finally have that $\omega$ is strictly increasing on $[0,\infty).$

Part 6) $\omega$ is sub-additive on $[0,\infty)$: We divide several cases to prove that \[\omega(s+t)\leq \omega(s)+\omega(t)~~{\rm  for}~ s,t\in [0,\infty).\]

Case I. If $s+t\geq 1$.  From  $\omega(\theta)\geq \theta$ for $\theta\in [0,\infty)$ we get \[\omega(s+t)=s+t\leq \omega(s)+\omega(t).\]

Case II. If $e^{-1}< s+t<1$.

Assume $s,t\in (e^{-1},1).$ Clearly $\omega(s+t)\leq \omega(s)+\omega(t).$

Assume $s\in (0,e^{-1}),t\in [e^{-1},1).$ Then
\begin{eqnarray*}
  &&\omega(s)+\omega(t)-\omega(s+t)\\&=& s(|\log s|^\alpha+1)+\frac{e-2}{e-1}t+\frac{1}{e-1}-\Big[\frac{e-2}{e-1}(s+t)+\frac{1}{e-1}\Big]\\
   &=& s\Big[|\log s|^\alpha+1-\frac{e-2}{e-1}\Big] \\
   &\geq &\Big[2-\frac{e-2}{e-1}\Big]s>0,
   \end{eqnarray*}
   and so $\omega(s+t)\leq \omega(s)+\omega(t).$

Assume $s, t\in (0,e^{-1}).$ Then
\begin{eqnarray*}
  &&\omega(s)+\omega(t)-\omega(s+t)\\
  &=& s(|\log s|^\alpha+1)+t(|\log t|^\alpha+1)-\left(\frac{e-2}{e-1}(s+t)+\frac{1}{e-1}\right)\\
   &=& s\Big[|\log s|^\alpha+1-\frac{e-2}{e-1}\Big]+t\Big[|\log t|^\alpha+1-\frac{e-2}{e-1}\Big]-\frac{1}{e-1}\\
   &\geq&(s+t)\frac{e}{e-1}-\frac{1}{e-1}> 0,
   \end{eqnarray*}
   and so $\omega(s+t)\leq \omega(s)+\omega(t).$

Case III. If $0<s+t\leq e^{-1}.$ Since \[f(t)=\frac{\omega(t)}{t}=\frac{t(|\log t|^\alpha+1)}{t}=|\log t|^\alpha+1\] is decreasing on $(0,e^{-1}]$ then
\[\frac{\omega(s+t)}{s+t}\leq \frac{\omega(s)}{s}, \frac{\omega(s+t)}{s+t}\leq \frac{\omega(t)}{t}.\]
Thus \[\omega(s+t)=\omega(s+t)\frac{s}{s+t}+\omega(s+t)\frac{t}{s+t}\leq \omega(s)+\omega(t).\]
\hfill$\square$

\begin{rem}

Here $\alpha=2$ is the best possible value for the monotonicity of this particular piecewise formula.
\end{rem}

\begin{ex}\label{pr1}
Let $X$ be a Banach space. Given $0<\alpha\leq 2$, let $\omega: [0,\infty)\to[0,\infty)$ be a function given by
\[\omega(t)=\begin{cases}0,&t=0\\
t(|\log t|^\alpha+1),&0<t<e^{-1}\\
\frac{e-2}{e-1}t+\frac{1}{e-1},&e^{-1}\leq t< 1\\
t,&t\geq 1.\\

\end{cases}\]
Then $\mathcal{F}_{\omega}(X)$ and ${\rm ker}\,{\beta_{\omega}}\oplus  X$  are $(\alpha,\alpha)$-almost Lipschitz homeomorphic.

\end{ex}
\pf
 First,  $\omega$ is a  nontrivial and strongly normalized gauge.  That is,  $\omega $ is a continuous increasing subadditive function with $\omega(0)=0$,  $\omega(t)\geq t$ for $0\leq t\leq {1}$ and  $\omega(t)=t$ for all $t\geq{1}$ and $\lim\limits_{t\to{0}}\frac{\omega(t)}{t}=\infty.$\par

By \cite[Proposition 5.1]{Kalton2004}, there is a uniform homeomorphism $\varphi: {\rm ker}\,{\beta_{\omega}}\oplus_{1} X\to\mathcal{F}_{\omega}(X)$ defined by
\[\varphi(\gamma, x)=\gamma+\delta_{\omega}(x)\]
such that $\varphi^{-1}(\gamma)=\left(\gamma-\delta_{\omega}\left(\beta_{\omega}(\gamma)\right),\beta_{\omega}(\gamma)\right)$ and
$\omega_{\varphi}(t)\leq2\omega(t), \,\omega_{\varphi^{-1}}(t)\leq3\omega(t),\, t>0.$
Thus $\mathcal{F}_{\omega}(X)$ and ${\rm ker}\,{\beta_{\omega}}\oplus_{1} X$  are $(\alpha,\alpha)$-almost Lipschitz homeomorphic.
\hfill$\square$
\par

In particular, by taking $X=c_0$ or $L_1[0,1]$ (or any separable, non-RNP Banach space) and considering the natural quotient map (the barycentric map) $\beta_{\omega}: \mathcal{F}_{\omega}(X)\to X$, we get the first pair of almost Lipschitz homeomorphic  but not Lipschitz homeomorphic
separable Banach spaces.
\begin{thm}\label{ex1}
Given $0<\alpha\leq 2$, let $\omega$ be as in Proposition \ref{PR1}. Then  separable Banach spaces $\mathcal{F}_{\omega}(c_0)$ and ${\rm ker}\,{\beta_{\omega}}\oplus_{1} c_0$ are $(\alpha,\alpha)$-almost Lipschitz homeomorphic but not Lipschitz homeomorphic.

\end{thm}
\pf
Recall that Aliaga, Gartland, Petitjean and Proch\'{a}zka \cite[after Theorem 4.6]{AGPP2022} showed that if $\omega$  is a nontrivial gauge then $\mathcal{F}_{\omega}(X)$ has the Radon-Nikod\'ym property. So $\mathcal{F}_{\omega}(c_0)$ has the Radon-Nikod\'ym property.  Note that the Radon-Nikod\'ym property is invariant under Lipschitz embeddings. Now $c_0$ is a separable Banach space without the Radon-Nikod\'ym
property,  it follows from the classical differentiability
results that $c_0$ cannot Lipschitz embed into $\mathcal{F}_{\omega}(c_0)$. Furthermore ${\rm ker}\,{\beta_{\omega}}\oplus_{1} c_0$ and $\mathcal{F}_{\omega}(c_0)$ are not Lipschitz homeomorphic.\hfill$\square$

\medskip
\noindent\textbf{Acknowledgments.} The authors thank Professor G. Lancien for helpful comments on an earlier version of the manuscript. They also thank the faculty and students of the Functional Analysis Seminar at Xiamen University for helpful discussions and suggestions. This work was supported by NSFC (Grant No. 12071389) and the Natural Science Foundation of Fujian Province (Grant No. 2024J01028).

\bibliographystyle{amsalpha}

\end{document}